\input amstex
\documentstyle{amsppt}
\document
\magnification=1200
\NoBlackBoxes
\nologo
\pageheight{18cm}


\bigskip

\centerline{\bf EXTENDED MODULAR OPERAD}

\medskip

\centerline{\bf A.~Losev${}^1$, Yu.~Manin${}^2$}

\medskip

\centerline{\it ${}^1$Institute of Theoretical and Experimental Physics, Moscow, Russia}

\smallskip

\centerline{\it ${}^2$Max--Planck--Institut f\"ur Mathematik, Bonn, Germany}

\smallskip

\centerline{\it and Northwestern University, Evanston, USA}

\bigskip

{\bf Abstract.} This paper is a sequel to [LoMa] where 
moduli spaces of painted stable curves were introduced and
studied. We define the extended modular operad of genus zero,
algebras over this operad, and study the formal differential geometric
structures related to these algebras: pencils of flat connections
and Frobenius manifolds without metric. We focus here
on the combinatorial aspects of the picture.
Algebraic geometric aspects are treated in [Ma2].

\bigskip

\centerline{\bf \S 0. Introduction and summary}

\medskip

This paper, together with [Ma2],  constitutes a sequel to [LoMa] where some
new moduli spaces of pointed curves were introduced and
studied. We start with a review of the main results
of [LoMa] and then give a summary of this paper.

\medskip

{\bf 0.1. Painted stable curves.} 
Let $S$ be a finite set. {\it A painting} of $S$
is a partition of $S$ into two disjoint subsets: {\it white} $W$
and {\it black} $B$.

\smallskip

Let $T$ be a scheme, $S$ a painted set, $g\ge 0$. 
{\it An $S$--pointed (or labeled) curve of genus $g$
over $T$} consists of the data
$$
(\pi:\,C\to T;\,x_i:\,T\to C,\ i\in S)
\eqno (0.1)
$$
where

\smallskip

(i) $\pi$ is a flat proper morphism whose geometric fibres $C_t$
are reduced and connected curves, with at most ordinary double
points as singularities, and $g=H^1(C_t,\Cal{O}_{C_t}).$

\smallskip

(ii) $x_i, i\in S,$ are sections of $\pi$ not containing
singular points of geometric fibres.

\smallskip

(iii) $x_i\cap x_j=\emptyset$ for all $i\in W,\ j\in S,\ i\neq j.$

\medskip

Such a curve $(0.1)$ is called {\it painted stable}, if the normalization
of any irreducible component $C^{\prime}$ of a geometric fibre
carries $\ge 3$ pairwise different special points when $C^{\prime}$ is of genus $0$
and $\ge 1$ special points when $C^{\prime}$ is of genus $1.$ 
Special points are inverse images of singular points
and of the structure sections $x_i$. Equivalently, such a normalization
has only a finite automorphism group fixing the special points.

\smallskip

We get the usual notion of a (family of) stable painted curves,
if $S$ consists of only white points.

\smallskip

The (dual) modular graph of a geometric fibre is defined
in the same way as in the usual case (for the conventions
we use see [Ma1], III.2). Tails now can be of two
types, we may refer to them and their marks 
as ``black'' and ``white'' ones as well, and call the graph
painted one. More on the geometry
of such graphs and morphisms between them see in the subsections 1.1--1.2
below. The isomorphism class of a painted modular
graph is called {\it the combinatorial type} of the respective curve.

\medskip

{\bf 0.2. Moduli stacks $\overline{L}_{g,S}.$} In sec. 4 of [LoMa], we
gave two equivalent descriptions of the stacks $\overline{L}_{g,S}$ where
$S$ is a painted set. In [Ma2], Lemma 1.2.2,  a third description was given
based upon a more general definition due to B. Hassett [H].
For reader's convenience, we will reproduce them here.

\medskip

(i) $\overline{L}_{g,S}$ consists of families of curves of {\it chain
stable combinatorial types.} 

\smallskip

To define the latter, consider first
combinatorial types of classical (semi)stable curves 
with only white points labeled by a finite set $W$. They
are isomorphism classes of graphs, whose vertices are
labeled by ``genera'' $g\ge 0$ and tails are bijectively labeled by elements
of $W$. Stability means that vertices of genus 0 bound $\ge 3$
flags, and vertices of genus 1 bound $\ge 1$ flags.
Graphs can have edges with only one vertex, that is, simple loops.

\smallskip

Starting with  the geometric
realization of such a graph $\Gamma$,  we can obtain an infinite series
of painted graphs, which will be called chain
stable. Namely,
subdivide edges and tails of $\Gamma$ by a finite set of new vertices
of genus zero (on each edge or tail, this set may be empty). If
a tail was subdivided, move the respective label (from $W$)
to the newly emerged tail. Distribute the black
tails labeled by elements of $B$ arbitrarily among the old and the new
vertices. The resulting graph is called  chain stable
if it becomes stable after repainting
black tails into white ones. The new vertices depict chains of
$\bold{P}^1$'s stabilized by black points. Each end of such a  chain
carries either a white tail or a singular special point (or both).

\medskip

(ii) The second description of $\overline{L}_{g,S}$ is based upon
the inductive construction
of ``adjoining a new black point''. 

\smallskip

Namely,
starting with an $S$--pointed family $(C/T, (x_i\,|\,i\in S))$ as in (0.1) 
we will produce
another $S'$--pointed family $(C'/T', (x_i^{\prime}\,|\,i\in S'))$
where $T'=C$, $S'=S\cup\{*\}$, $*$ being a new black label.
The new curve and sections
will be produced in two steps. At the first step
we make in (0.1) the base change $C\to T$.
We get an $S$--pointed curve $X:=C\times_S C$,
with sections $x_{i,C}.$
We then add an extra section $\Delta :\,C\to C\times_S C$
which is the relative diagonal, and mark it by $*$.
We have not yet produced an $S^{\prime}$--pointed curve over
$T^{\prime}=C$, because the extra black section  will
intersect singular points of the fibres and
white sections as well. 

\smallskip

Now comes the second step of the construction, 
where we birationally
modify $C\times_T C\to C$ as in [Kn], Definition 2.3.  
Define $C^{\prime}:= \roman{Proj\,Sym}\,\Cal{K}$
as the relative projective spectrum of the symmetric algebra
of the sheaf $\Cal{K}$ on $X=C\times_T C$ which is defined as
the cokernel of the map
$$
\delta :\,\Cal{O}_X\to  \Cal{J}_{\Delta}\check{}\,\oplus\,
\Cal{O}_X(\sum_{i\in A}x_{i,C}),\ \delta (f) = (f,f).
\eqno(0.3)
$$
Here $\Cal{J}_{\Delta}$ is the $\Cal{O}_X$--ideal of $\Delta$,
and $\Cal{J}_{\Delta}\check{}$ is its dual sheaf
considered as a subsheaf of meromorphic
functions, as in [Kn], Lemma 2.2 and Appendix.

\smallskip

As a result, we get an $S^{\prime}$--pointed curve,
because Knudsen's treatment of his modification is local and
can be directly extended to our case. If the initial curve
was painted stable, the new one will be painted stable as well.

\smallskip

Moreover, this construction is obviously functorial, in the
sense that it produces a stack if we started with a stack.

\smallskip

We are now ready to construct $\overline{L}_{g,S}.$ For concreteness,
we assume that $W=\{1,\dots ,m\}$ and $B=\{1,\dots ,n\}$ are initial segments
of (two copies of) natural numbers.

\smallskip

If $g\ge 2$, $m\ge 0$, we start
with $\overline{M}_{g;m}=\overline{M}_{g;m,\emptyset}$ and
adjoin $n$ new black points, one in turn. Denote the resulting
stack by $\overline{L}_{g;m,n}.$ 

\smallskip

For $g=1$, we repeat this construction for $m\ge 1.$  
We need one more sequence of stacks, corresponding to $m=0$
(elliptic curves stabilized only with black points).  
Since we want to restrict ourselves to Deligne--Mumford stacks, we 
start with $\overline{M}_{1;0,1}$ which is defined by repainting
the white point in $\overline{M}_{1;1}$, and adjoin
black points to get the sequence
$\overline{L}_{1;0,n}$, $n\ge 1.$ 

\smallskip

Finally, for $g=0$ we obtain the series of spaces $\overline{L}_n=\overline{L}_{0;2,n}$,
thoroughly studied in [LoMa]
and moreover $\overline{L}_{0;m,n}$, for all $m\ge 3, n\ge 0$.

\smallskip

From our construction it follows that the genus zero stacks are actually
smooth manifolds.

\medskip

(iii) The third, and technically most useful, description identifies
$\overline{L}_{g,S}$ with certain stacks of {\it weighted pointed stable curves}
studied in [H], see also [Ma2], Definition 1.2.1. Roughly speaking,
in Hassett's definition we endow each section $x_s$ with a weight 
$a_s\in \bold{Q}\,\cap\,(0,1].$ {\it Weighted stability} means that 
(an appropriate power of) $\omega_{C/T} (\sum_s a_sx_s)$ is relatively ample,
and moreover, the sum of weights of any subfamily of $x_s$ coinciding
on a geometric fiber does not exceed 1. Lemma 1.2.1 in [Ma2] shows that
if $S$ is painted, and if we choose weights so that $a_s=1$ for $s\in W$
and $\sum_{t\in B} a_t \le 1$, the two definitions of stability
coincide.

\medskip

{\bf 0.3. Extended modular operad.} The classical modular spaces
$\overline{M}_{g,m+1}$ come together with the action of $\bold{S}_{m+1}$. Moreover, 
there is a family of clutching morphisms
$$
\overline{M}_{g_1,m_1+1}\times \overline{M}_{g_2,m_2+1}\to
\overline{M}_{g_1+g_2,m_1+m_2+1}
\eqno(0.4)
$$
which on the level of geometric points can be described
as gluing one (say, the last) labeled point
of the first pointed curve to, say, the first labeled point
of the second curve. Finally, there is a similar
family of morphisms 
$$
\overline{M}_{g-1,m+1} \to \overline{M}_{g,m-1}
\eqno(0.5)
$$
gluing together two different labeled points.

\smallskip

We can  axiomatize
the evident relations between these structures and get
the abstract notion of the modular operad in the category
of $DM$--stacks. One can also add a linear structure,
by passing to the (homology) motives or any of their
realizations and get the notion of the modular operad in
a tensor category, as well as the basic example
of this notion. This was done by E.~Getzler and M.~Kapranov 
in [GeK2]. Restricting ourselves to the $g=0$ case, we
get almost a classical operad, with two modifications:
first, its $m$--th component $\overline{M}_{g,m+1}$
is acted upon by $\bold{S}_{m+1}$ rather than by $\bold{S}_{m}$;
second, it lacks the $m=1$--component since $\bold{P}^1$
with only two labeled points is unstable. 

\smallskip

The first modification means that we actually deal with {\it a cyclic
operad}, the notion which was  as well introduced and studied by
E.~Getzler and M.~Kapranov ([GeK1]). Recall that generally
classical linear operads classify $k$--algebras of a given type,
like associative, commutative, Lie, Poisson etc. The component $P(m)$
of an operad $P$ (in our case $H_*(\overline{M}_{0,m+1})$ )
is the linear space of polylinear
operations of $m$ arguments that can be constructed from the
basic operations (multiplications, commutators etc),
modulo identical relations supplied by the axioms
(associativity, Jacobi etc). The action of $\bold{S}_m$ upon $P(m)$
corresponds to the permutations of the arguments.
If $A$ is an algebra over $P$, each $p\in P(m)$ defines thus
a tensor $p:\,A^{\otimes m}\to A,$ and we know how these tensors
transform under all permutations of $m$ lower indices.
If $A$ is endowed with a scalar product, we can lower
the $m+1$--th index as well identifying
$p$ with a map $A^{m+1}\to k$. It makes sense to require  
the scalar product and the action of $\bold{S}_{m+1}$ on such tensors
to be a part of the structure of algebras
of a given type. Cyclic operads then are well--suited
for axiomatizing such a situation.

\smallskip

We now pass to the problem of the 1--component. It is missing
because a projective line with two labeled (white)
points has an infinite automorphism group $G_m$ and so
is unstable. The trick suggested in [LoMa]
consists in stabilizing $\bold{P}^1$ with an arbitrary number
of additional black points. This forces us to add
degenerations  which are arbitrary chains of such $\bold{P}^1$'s.
Thus our candidate
for the missing $\overline{M}_{0,2}$ is
$\coprod_{n\ge 1} \overline{L}_{0,2,n}$. The remaining
spaces $\overline{L}_{g,m,n}$ constitute a minimal extension 
of the family $\overline{M}_{g,m}$ containing all
$\overline{L}_{0,2,n}$ and stable with respect to  the
operadic (clutching) morphisms.

\smallskip

Specifically, we put
$$
\widetilde{L}_{g,m}:=\coprod_{n} \overline{L}_{g,m+1,n} \, .
\eqno(0.6)
$$
Morphisms (0.4) and (0.5) in which we allow to glue together {\it only white points}
induce the clutching maps
$$
\widetilde{L}_{g_1,m_1}\times \widetilde{L}_{g_2,m_2}\to
\widetilde{L}_{g_1+g_2,m_1+m_2} \, ,
\eqno(0.7)
$$
$$
\widetilde{L}_{g-1,m+1} \to \widetilde{L}_{g,m-1}\, .
\eqno(0.8)
$$
The symmetry group $\bold{S}_{m+1}$ acts upon $\widetilde{L}_{g,m}$
by renumbering white points. 
\smallskip

Moreover, each stack $\widetilde{L}_{g,m}$ is
a union of infinitely many irreducible $DM$--stacks which are indexed
by the number of black points; the group $\bold{S}_n$
acts on the $n$--th component by renumbering black points.
These structures induce the respective symmetries
of the linearized versions. This is the price we paid for
acquiring the $m=1$ component. It has a nice meaning in the language
of the respective algebras: an algebra now consists of
two spaces $(A,T)$; polylinear operations are of the form
$A^{\otimes m}\otimes T^{\otimes n}\to A$. They
depend on $t_i\in T$ as parameters which are associated with black
points. Thus we acquired degrees of freedom
allowing us to encompass certain deformations. 
\smallskip

The key result of the theory of the (cyclic) homology operad $\{ H_*(\overline{M}_{0,m+1})\}$
consists in establishing  a bijection between finite--dimensional
operadic algebras and formal Frobenius manifolds.
A more detailed analysis
which we relegate to a future paper shows that in our present context we get
exactly the deformations associated with gravitational
descendants.

\medskip

{\bf 0.4. Genus zero extended operad.} It is convenient to imagine
its bigraded components $\overline{L}_{0,m+1,n}$ sitting at the lattice points
$(m,n)$ of the first quadrant. Then the lowest horizontal
line is the classical genus zero modular operad $\overline{M}_{0,m+1}$
whereas the leftmost vertical line represents the
1--component of the extended operad. Operadic multiplication
must induce a monoid structure on the 1--component, and in our case
it is supplied by the appropriate clutching morphisms described above and 
more formally in [Ma2], [LoMa]:
$$
\widetilde{L}_{0,m_1+1}\times \widetilde{L}_{0,m_2+1}\to
\widetilde{L}_{0,m_1+m_2+1} \, ,
\eqno(0.9)
$$

\smallskip

In [LoMa], it was proved that the appropriately defined
algebras over the 1--component $\widetilde{L}_{0,1}$ bijectively
correspond
to the  pencils of formal flat connections. The treatment in [LoMa]
stressed also the toric picture of $\widetilde{L}_{0,1}$;
for certain generalizations of this picture, cf. [R].

\smallskip

Thus at the boundary of the $(m,n)$--quadrant we have
a pretty detailed understanding of the geometry of the respective
moduli spaces, algebraic structure of their homology and
cohomology including the operadic formalism, and finally,
a reduction of the theory of operadic algebras to the
study of specific differential equations.

\smallskip

The main goal of [Ma2] and this paper taken together is to extend these results from
the boundary of the $(m,n)$--quadrant to the whole quadrant.

\smallskip

Here we focus mainly on the combinatorial aspects of the picture;
those arguments for which algebraic geometry is indispensable,
are treated in [Ma2].

\medskip

{\bf 0.5. Contents of the paper.} In \S 1 and  \S 2 we define for
any painted set $S$
purely combinatorially a graded ring $H^*_S$ and a graded module
$H_{*S}$ over this ring and study their structure. Using
the results of this study and additional geometric
arguments, we establish in [Ma2] canonical isomorphisms
$H^*_S=H^*(\overline{L}_{0S})$ and $H_{*S}=H_*(\overline{L}_{0S})$,
together with identification of the action of $H^*_S$
upon $H_{*S}$ as the cap product.

\smallskip

In \S 3 we extend the construction of $H^*$ and $H_*$ to the labeled trees
replacing $S$. Algebraic geometric meaning of this
extension consists in studying (co)homology
of moduli spaces of painted stable curves of a given
combinatorial type and their degenerations. Combinatorially,
these rings/modules together with various
morphisms connecting them constitute a version
of graded (co)operad.

\smallskip

In \S 4 we introduce several versions of algebras over this (co)operad
and describe them in terms of correlators and top correlators.

\smallskip

Finally, in \S 5 we demonstrate the  geometric meaning
of the differential equations satisfied by the generating
series for top correlators. There are two basic geometric languages,
corresponding to what we call in \S 5  {\it the commutativity}
and {\it the oriented associativity equations}. One language leads to the study
of pencils of flat connections as in [LoMa1],
another to a version of ``Frobenius manifolds with affine flat
structure but without invariant metric''. Both versions are related
by the induction of flat affine structure from a
fiber. This procedure first appeared in [Lo1],[Lo2].

\smallskip

{\it Acknoledgements.} The research of A.~L. was partially supported by
grants  RFBR 01-01-00548,
by INTAS-99-590 and by support  of scientific schools grant 00-15-96557.

\bigskip

\centerline{\bf \S 1. Combinatorial cohomology rings}

\medskip

{\bf 1.1. Graphs.} We define {\it a graph $\tau$}
as a quadruple $(V_{\tau}, F_{\tau}, \partial_{\tau},j_{\tau})$
where $V_{\tau}$, resp. $F_{\tau}$, are finite sets of
vertices, resp. flags; $\partial_{\tau}:\, F_{\tau}\to V_{\tau}$
is the boundary, or incidence, map; $j_{\tau}:\,F_{\tau}\to F_{\tau}$ is an involution of the set of flags. 
This combinatorial definition is related to the more common one via
the notion of
{\it the geometric realization
of $\tau$}. The latter  is a topological space which is obtained from
$F_{\tau}$ copies of $[0,1]$ by gluing together
points $0$ in the copies corresponding to each vertex $v\in V_{\tau}$,
and by gluing together points $1$ in each orbit of $j_{\tau}$.
This motivates the introduction of the following auxiliary sets
and their geometric realizations: the set $E_{\tau}$ of {\it edges} of
$\tau$, formally consisting of cardinality two orbits of $j_{\tau}$,
and the set $T_{\tau}$ of {\it tails}, consisting of those
flags, which are $j_{\tau}$--invariant. 

\smallskip

We will very often think and speak about graphs directly in terms
of their geometric realizations. In particular, $\tau$ will be called
{\it connected} (resp. {\it tree}, resp. {\it forest})
if its geometric realization is connected (resp. connected and
has no loops, resp. is 
a disjoint union of trees).

\smallskip

Situations in which one must be more careful and invoke
the initial combinatorial definition  usually
occur when one has to consider morphisms of graphs.
We will define here three types of generating morphisms
which we will need. For more sophisticated definitions,
cf. [BeMa]. 

\smallskip

In fact, [BeMa] considers graphs endowed with additional labeling
of their vertices with integers $\ge 0.$ Such graphs arise
as dual graphs of (semi)stable pointed curves: vertices correspond to the irreducible
components, labels to the genera of their normalizations.
This part of the structure can be neglected in most of this paper
which treats only the combinatorics of the genus zero
curves whose dual graphs are trees with all vertices labeled by zero.

\smallskip

Moreover, flags are in a bijection with the special points of the
normalized curve. Hence when $S$ is weighted and/or painted,
we must consider graphs whose flags are weighted, resp. painted.

\medskip

(i) {\it An isomorphism} $\sigma\to\tau$ consists of two
bijections  $\varphi^F:\,F_{\tau}\to F_{\sigma}$, 
$\varphi_V:\,V_{\sigma}\to V_{\tau}.$
They must be compatible with $\partial_{\sigma}$,
$\partial_{\tau}$,  and $j_{\sigma}$, $j_{\tau}$. 
As a corollary, they induce bijections of edges and  tails as well.
(Notice that maps between flags go in reverse
direction with respect to the morphisms of graphs.
This  convention will persist below.)

\smallskip

(ii) Let $\sigma$ be a graph, and $e$ its edge
consisting of two flags $e_1,e_2.$
A morphism $c:\,\sigma\to\tau$ {\it contracting $e$}
geometrically simply deletes the interior of $e$ if it is a loop
with vertex $v$,
or else contracts $e$ so that its ends $v_1,v_2$
become a new vertex $v$ otherwise.
Combinatorially, it is given by a map $c^F:\,F_{\tau}\to F_{\sigma}$
which identifies $F_{\tau}$ with $F_{\sigma}\setminus \{e_1,e_2\}$
and is compatible with involutions, and a map $c_V:\,V_{\sigma}\to V_{\tau}$
which glues together vertices of $e_1,e_2$ and is a bijection
elsewhere.
Distribution of flags between all vertices
except for those of $e_1,e_2$ must be compatible
with these maps. Flags of $\tau$ at $v$ (if $e$ is not a loop) 
are thus $F_{\sigma}(v_1)\cup F_{\sigma}(v_2)\setminus \{e_1,e_2\}.$

\smallskip

(iii) Let $\sigma$ be a graph and $\{f_1,f_2\}$ its two
different tails. A morphism $\sigma\to \tau$ which {\it glues} these
two tails together consists of two bijections $F_{\tau}\to  F_{\sigma}$, $V_{\sigma}\to V_{\tau}$ compatible with $\partial_{\sigma}$,
$\partial_{\tau}$, and making of $\{f_1,f_2\}$
one new edge.

\smallskip

Composition of morphisms is the composition of the relevant maps
(in reverse order for flags). A composition of contractions of edges and isomorphisms will be called simply a contraction. 
The category of trees and contractions will reappear in \S 3.

\medskip

{\bf 1.1.1. Painted graphs and labeled graphs.} {\it A painting}
of $\tau$ is a partition of its tails into two subsets:
white and black. If a painting is not given
explicitly, we assume that all tails are white.

\smallskip

Morphisms between painted graphs are subject to the following
restrictions: isomorphisms and contractions do not change color of tails;
it is allowed to glue together only white tails.

\smallskip

Let $S$ be a set, and $\tau$ a graph. {\it An $S$--marking}
(or labeling) of $\tau$ is a bijection $S\to T_{\tau}.$
Equivalently, we call a graph with an $S$--marking
an $S$--graph, and a morphism between two $S$--graphs
identical on $S$, an $S$--morphism.
We identify a painting of $\tau$ with a painting of $S$.

\medskip

{\bf 1.1.2. Stable and stably painted trees.} Let now $\tau$ be a tree.
We recall that it is called {\it stable}, iff its
every vertex is incident to at least three flags.

\smallskip

A painted tree is called {\it painted stable}, iff

\smallskip

(i) It is stable.
\smallskip

(ii) Every end vertex of $\tau$ carries at least one white tail.
(An end vertex is a vertex incident to 1 or 0 edges).

\medskip

In particular, an one edge $S$--tree is the same as an
unordered 2--partition $S=\sigma_1\cup \sigma_2$
showing how the tails are distributed between the two vertices.
Painted stability means that $\sigma_1$ and $\sigma_2$ each contain
at least one white point, and $|\sigma_i|\ge 2.$

\smallskip

We will often use the following notation. If $S_1,S_2$ are two subsets
of $S$ and $\sigma$ is a 2--partition of $S$, then
$S_1\sigma S_2$ is a shorthand for the following
statement: all elements of $S_1$ lie in one part of
$\sigma$ whereas all elements of $S_2$ lie in another one.
If one of the sets is empty, we may omit it.

\medskip

{\bf 1.2. 2--partitions and trees.}
Let $\tau$ be a tree, $e$ its edge. If we delete the interior
of (the geometric realization of)
$e$, the tree will split into two connected components,
and the set of its tails will break into two subsets.
If $\tau$ is $S$--marked, we get thus an unordered 2--partition $\sigma_e$ of $S$.
If $\tau$ is painted stable, this partition will also be
painted stable.

\smallskip

Let $\sigma =(\sigma_1,\sigma_2)$ and $\tau =(\tau_1,\tau_2)$ 
be two non--trivial unordered partitions of $S$. Denote
by $\delta (\sigma ,\tau )$  the number of non--empty
pairwise distinct sets among $\sigma_k\cap\tau_l$
diminished by 2. 

\smallskip

Clearly, $\delta (\sigma ,\tau )$ can take only values
$0,1,2.$ Moreover,

\smallskip

(i) $\delta (\sigma ,\tau )=0$ {\it iff} $\sigma =\tau .$

\smallskip

(ii) $\delta (\sigma ,\tau )=1$ {\it iff the following condition is satisfied.
There exists a partition of $S$ into three
(non-empty pairwise disjoint) subsets $\rho_i,\,i=1,2,3,$
such that one of the partitions $\sigma,\tau$ is $(\rho_1\cup\rho_2,\rho_3)$,
and another is $(\rho_1,\rho_2\cup\rho_3).$ The partition
$(\rho_i)$ is determined uniquely up to a permutation of $\rho_1$ and $\rho_3$.}

\smallskip

To see this, notice that there must be exactly one empty
intersection $\sigma_a\cap\tau_b$. Renumbering the sets if necessary
we may assume that $\sigma_2\cap\tau_1=\emptyset$. Then 
$\rho_1:=\tau_1, \rho_2:=\sigma_1\cap\tau_2$ and $\rho_3:=
\sigma_2$ will do.

\smallskip

Moreover, assume in this situation that $\sigma$ and $\tau$
are painted stable as one--edge trees. Denote by
$\sigma *\tau$ the two--edge $S$--tree whose tails are distributed
between its three vertices according to the partition
$(\rho_i)$, with $\rho_2$ in the middle. Then
$\sigma *\tau$ is stable painted as well. In fact,
its end vertices carry at least two tails of which at least one is white because
this holds for $\sigma$ and $\tau$; and its middle vertex
carries at least three flags because $\rho_2$ is non--empty.

\smallskip

(iii) {\it   $\delta (\sigma ,\tau )=2$ iff there exists
an ordered quadruple of pairwise distinct elements
$i,j,k,l\in S$ such that $ji\sigma kl$ and $jk\tau il.$
Moreover, if $S$ is painted and $\sigma$, $\tau$ are painted
stable, $i,j,k,l$ can be chosen in such a way that
the partitions induced by $\sigma$ and $\tau$
on this quadruple are also painted stable.}

\smallskip

In fact, if one disregards stability, it suffices to choose $j\in \sigma_1\cap\tau_1,$
$i\in \sigma_1\cap\tau_2,$ $k\in \sigma_2\cap\tau_1,$ 
$l\in \sigma_2\cap\tau_2.$ This is possible
exactly when $\delta (\sigma ,\tau )=2$.

\smallskip

If $\sigma$ and $\tau$ are painted stable,
consider the $(2,2)$ matrix of sets $(\sigma_a\cap\tau_b).$
A contemplation will convince the reader that
at least one of the two diagonals of this matrix
has the following property: each of its elements contains
a white label. Hence either $i$ and $k$, or else
$j$ and $l$ can be chosen white. In both cases,
the induced partitions of $\{i,j,k,l\}$
are painted stable.

\smallskip

We will call a set $E$ of pairwise distinct (painted) stable 2--partitions 
$\{\sigma\}$ of a (painted) stable set $S$
{\it good}, if $\delta (\sigma,\sigma')=1$ for any 
$\sigma\ne\sigma'$ in $E$. In particular, empty set
is good, and any one--element set is good.
 
\medskip

\proclaim{\quad 1.2.1. Lemma} (a) For any (painted) stable
$S$--tree $\tau$, the set $\{\sigma_e\,|\,e\in E_{\tau}\}$ is good.

\smallskip

(b) The map $\tau\mapsto \{\sigma_e\,|\,e\in E_{\tau}\}$ establishes
a bijection between (painted) stable $S$--trees 
up to an $S$--isomorphism, and good families.

\smallskip

(c) Let $E'$ be a non--empty good family, $\sigma\in E'$,
$E=E'\setminus \{\sigma\}$. Let $\tau'$, resp. $\tau$,
be the respective $S$--trees, and $e\in E_{\tau^{\prime}}$ the edge
corresponding to $\sigma$. Then there is a unique
$S$--morphism $\tau^{\prime}\to\tau$ contracting $e$
and inducing the tautological embedding $E\subset E'.$
\endproclaim

\smallskip

{\bf Proof.} If one disregards painting, statements (a) and (b)
are proved
in the Proposition 3.5.2  of Chapter III of [Ma1]. To deal
with painted trees,
it remains only to check that if one starts with a non--empty good
family and produces the respective tree, it will be
painted stable, that is, its every end vertex carries a 
white tail. In fact, such a vertex belongs to a unique
edge $e$, and the tails carried by it constitute
a part of the partition $\sigma_e$. Since $\sigma_e$
is painted stable, the conclusion follows.

\smallskip

We leave (c) as an exercise to the reader.

\medskip

\proclaim{\quad 1.2.2. Proposition} Let $\tau$ be a painted stable
$S$--tree, and $\sigma$ a painted stable 2--partition of $S$.
Then exactly one of the following alternatives takes place.

\smallskip

(i) There exists an edge $e\in E_{\tau}$ such that $\sigma =\sigma_e.$
In this case we will say that $\sigma$ breaks $\tau$ at $e$.

\smallskip

(ii) There exists a vertex $v\in V_{\tau}$, an $S$--tree $\tau^{\prime}$
and an $S$--morphism $\tau^{\prime}\to\tau$ contracting
an edge $e$ of $\tau'$ to the vertex $v$ of $\tau$ such that $\sigma =\sigma_e$.
In this case we will say that $\sigma$ breaks $\tau$ at 
$v.$

\smallskip

(iii) None of the above. In this case there is an edge
$e\in  E_{\tau}$ such that $\delta (\sigma ,\sigma_e)=2.$
We will say that $\sigma$ does not break $\tau$.
\endproclaim

\smallskip

{\bf Proof.} Consider the family of numbers $\delta (\sigma ,\sigma_e),$
$e\in E_{\tau}.$ If it contains zero, we are in the case (i), and then
all other numbers are 1, because of Lemma 1.2.1 (a).
If it consists only of 1's, we can add $\sigma$ to $\{\sigma_e\}$
and get a good family, which produces a new tree $\tau'$
and a contracting morphism $\tau'\to\tau$ with
required properties. If it contains 2, we are in the case (iii).

\medskip

{\bf 1.2.3. Remark.} In the following it will be useful
to keep in mind a more detailed picture of various
possibilities that can arise in the case (ii).

\smallskip

If we delete a vertex $v$ from the geometric realization of 
$\tau$, it will break into a set of $\ge 3$
connected components which we will call
{\it branches} of $\tau$ at $v$. Their set
is canonically bijective to the set
of flags $F_{\tau}(v)$ incident to $v$: we can say that
the  branch starts with the respective flag. 
In the extreme case, a branch can be a single tail.

\smallskip

Each branch
$b$ carries its full subset of tails, or their labels, $S_b$,
so that all $S_b$ together form a partition of $S$. A 2--partition $\sigma$
of $S$ breaks $\tau$ at $v$ precisely when each of its parts
is a union of several $S_b$'s. In this way
we establish a bijection between such 2--partitions
$\sigma$ of $S$ and 2--partitions $\alpha$ of $F_{\tau}(v).$
We will say that $\alpha$ is (painted) stable if the respective
$\sigma$ is.
Most of $\alpha$'s will be automatically
painted stable. Unstable $\alpha$'s are of two types:
\smallskip

(i) one part of $\alpha$ consists of one white tail; or

\smallskip

(ii) one part of $\alpha$ consists of only black tails.

\smallskip

In fact, if a branch starts with an edge, then
at an end vertex of this edge there are at least two tails,
of which at least one is white.

\medskip

{\bf 1.3. Combinatorial cohomology.} In the following
we fix a painted set $S$ with $|S|\ge 3$ containing at least two
white elements, and a commutative coefficient ring $k$.
Consider the family of independent commuting variables
$\{l_{\sigma}\}$ indexed by painted stable unordered
2--partitions $\sigma$ of $S$ and put $\Cal{R}_S:=k[l_{\sigma}].$
\smallskip

Call an ordered quadruple of pairwise distinct
elements $i,j,k,l\in S$ {\it allowed}, if
both partitions $ij|kl$ and $kj|il$ are painted stable.
For a painted stable $\sigma$ put $\epsilon (\sigma ;i,j,k,l)=1$
if $\{i,j,k,l\}$ is allowed and $ij\sigma kl$; $-1$,
if $\{i,j,k,l\}$ is allowed and $kj\sigma il$;
and $0$ otherwise.

\smallskip

Define
$$
R_{ijkl}:=\sum_{\sigma} \epsilon (\sigma ;i,j,k,l)\,l_{\sigma} \in \Cal{R}_S ,
\eqno(1.1)
$$
$$
R_{\sigma\sigma'}:=l_{\sigma}l_{\sigma'}\ if\ there\ exists\
an\ allowed\ quadruple\ i,j,k,l\ with\ ij\sigma kl,\,kj\sigma' il .
\eqno(1.2)
$$
Here $\sigma,\sigma'$ run over 2--partitions.
Denote by $I_S\subset \Cal{R}_S$ the ideal generated by all
elements (1.1) and (1.2) and define {\it the combinatorial cohomology
ring} by
$$
H_S^*:=\Cal{R}_S/I_S.
\eqno(1.3)
$$

\smallskip

The main result of [Ma2], Theorem 3.1.1, establishes a canonical isomorphism
of $H_S^*$ with the cohomology and the Chow ring of $\overline{L}_{0,S}.$
The proof heavily uses calculations of this paper and
the Keel's theorem who proved this result in [Ke] for the basic case when there are
no black points.

\smallskip

Returning to our combinatorial setup, let now $\tau$ be a painted stable $S$--tree. Put
$$
m(\tau ):=\prod_{e\in E_{\tau}}l_{\sigma_e} \in \Cal{R}_S
\eqno(1.4)
$$
and call such a monomial {\it good}. It depends
only on the $S$--isomorphism class of $\tau$. If $\tau$ is one--vertex
tree, we put $m(\tau )=1.$

\medskip

{\bf 1.3.1. Multiplication table.} Below we will calculate 
$l_{\sigma} m(\tau )\,\roman{mod}\, I_S.$ The answer will depend on 
the mutual position of $\sigma$ and $\tau$ as described in Prop. 1.2.2
which we will consider in reverse order.

\smallskip

{\it Case (iii): $\sigma$ does not break $\tau$.} Then
$$
l_{\sigma}m(\tau )\equiv 0\,\roman{mod}\,I_S
\eqno(1.5)
$$
because in this case $l_{\sigma}m(\tau )$ is divisible
by one of the elements (1.2).

\smallskip

{\it Case (ii): $\sigma$  breaks $\tau$ at a vertex $v$.} Then
$$
l_{\sigma}m(\tau )=m(\sigma*\tau)
\eqno(1.6)
$$
where $\sigma*\tau$ is the painted stable $S$--tree
corresponding to the set of 2--partitions $\{\sigma_e\,|\,e\in E_{\tau}\}
\,\cup\,\{\sigma\}.$

\smallskip

{\it Case (i): $\sigma$  breaks $\tau$ at an edge $e$.} In this case we will write several formulas depending on additional choices.
Namely, denote by $v_1,v_2$ two vertices of $e$ and by
$e_1$, $e_2$ the respective flags (``halves of $e$'').
Choose two distinct flags $I,J\in F_{\tau}(v_1)\setminus \{e_1\}$ and 
similarly $K,L\in F_{\tau}(v_2)\setminus \{e_2\}$ with the following
property:

\smallskip

(*) {\it there exist labels $i,j,k,l$ on the branches
starting with $I,J,K,L$ respectively,
forming an allowed quadruple.}

\smallskip

Again, we will call such $I,J,K,L$ an allowed quadruple.
It is important to notice that allowed $I,J,K,L$ always exist:
it suffices to take for $J$ any branch at $v_1$
carrying a white tail, for $L$ any branch at $v_2$
carrying a white tail, and complement them by arbitrary $I,K$.

\smallskip

Now we state that for each allowed quadruple of flags as above we have:
$$
l_{\sigma}m(\tau ) \equiv
-\sum_{\alpha :\,IJ\alpha e_1} m(\tau (\alpha ))-
\sum_{\beta :\,e_2\beta KL} m(\tau (\beta))\ \roman{mod}\,I_S.
\eqno(1.7)
$$
Here $\alpha$ (resp. $\beta$) runs over painted stable
2--partitions
of $F_{\tau}(v_1)$ (resp. $F_{\tau}(v_2)$) in the sense
of 1.2.3. If $\alpha$
determines a 2--partition $\sigma$ of $S$,
we denote by $\tau(\alpha )$ the former $\sigma*\tau$. 

\smallskip

To prove (1.7), we choose $i,j,k,l$ as in (*) and write first of all
$$
R_{ijkl}m(\tau )=\sum_{\rho} \epsilon (\rho ;i,j,k,l)\,l_{\rho} m(\tau ) \in
I_S.
\eqno(1.8)
$$
where we sum over 2--partitions $\rho$ of $S$. Now, if
$\epsilon (\rho ;i,j,k,l)= -1,$ then $kj\rho il$,
and since $ij\sigma kl$, we have $l_{\rho}l_{\sigma}\in I_S$;
but $l_{\sigma}$ divides $m(\tau )$, and so $l_{\sigma}m(\tau )\in I_S.$
Hence in the r.h.s. of (1.8) we may sum only over
$\rho$'s with $ij\rho kl$, and $\sigma$ is among this set.
Therefore,
$$
l_{\sigma} m(\tau )\equiv - \sum \Sb\rho :\,ij\rho kl \\ \rho\ne \sigma \endSb
l_{\rho} m(\tau )\,\roman{mod}\,I_S.
\eqno(1.9)
$$
In the sum (1.9) we may and will omit $\rho$'s that do not break
$\tau$. Moreover, no $\rho$ can break $\tau$ at an edge,
or at a vertex distinct from $v_1$ or $v_2$: otherwise
the condition $ij\rho kl$ cannot hold. Finally,
$\rho$'s which break $\tau$ at $v_1$ or $v_2$ will produce
precisely the right hand side of (1.7).

\medskip

{\bf 1.3.2. Corollary.} {\it The ring $\Cal{R}_S$ is linearly
spanned by $I_S$ and good monomials $m(\tau )$. The cohomology
ring $H^*_S$ is spanned by the classes of $m(\tau ).$}

\smallskip

In fact, the span of $I_S$ and $m(\tau )$ contains $1$
and is stable with respect to multiplications by all $l_{\sigma}.$

\smallskip

However, the space of linear relations
$I_S\,\cap\,\langle m(\tau )\rangle$ is generally
non--trivial. Below we will exhibit some elements of this
space which we will call {\it the standard ones}
and which will later be shown to span all the relations.

\medskip

{\bf 1.4. Standard relations between good monomials.} 
As above, we will calculate some expressions $R_{ijkl}m(\tau )$.
However,  formerly we started with a tree $\tau$, an edge
$e=e_{\sigma}$ and flags $I,J$ and $K,L$ at two ends of this edge.

\smallskip

Now we will start with a tree $\tau$, a vertex $v$,
an allowed quadruple of flags $I,J,K,L$ incident
to this vertex, and an allowed
quadruple of labels $i,j,k,l$ carried
by the respective branches. Looking at the summands at
the right hand side of (1.8),
we see that if $\rho$ does not break $\tau$, we can
omit it, and if $\rho$ breaks $\tau$ at an edge or at
a vertex different from $v$, then
$\epsilon (\rho ;i,j,k,l)=0.$ 

\smallskip

Partitions $\rho$ breaking $\tau$ at $v$ correspond
to 2--partitions $\alpha$  of $F_{\tau}(v).$ So finally we get
the standard relation
$$
\sum_{\alpha} \epsilon (\alpha ;I,J,K,L)\, m(\tau (\alpha )) \equiv 0\,
\roman{mod}\, I_S
\eqno(1.10)
$$
where the sign $\epsilon$ is $\pm 1$ depending on whether
$IJ\alpha KL$ or $KJ\alpha IL$.

\bigskip

\centerline{\bf \S 2. Combinatorial homology modules}

\medskip

{\bf 2.1. Homology space.} We keep notations of \S 1,
in particular, fix a painted set $S$ and a coefficient
ring $k$. Define a free $k$--module $\Cal{M}_S:=
\oplus_{\tau} k\,\mu(\tau )$ freely
generated by all $S$--isomorphism classes of painted stable
$S$--trees $\tau$. Denote by $J_S\subset\Cal{M}_S$
the submodule spanned by the standard relations (1.10)
in which $m(\tau)$ are replaced by $\mu (\tau )$:
$$
r(\tau ,v;I,J,K,L):= \sum_{\alpha} \epsilon (\alpha ;I,J,K,L)\, \mu (\tau (\alpha )) \in
J_S.
\eqno(2.1)
$$
Finally, define the homology $k$--space as
$$
H_{*S}:= \Cal{M}_S/J_S
\eqno(2.2)
$$
and put $[\mu (\tau )]:=\mu (\tau )\,\roman{mod}\,J_S$.
\smallskip

The main content of this section is the definition of an action
of $H^*_S$ on $H_{*S}$ which turns the homology
module
into a free module of rank one over the cohomology ring.
To be more precise, we start with defining an action of
the ring generators $l_{\sigma}$ upon $\Cal{M}_S$ which imitates the multiplication
table (1.5)--(1.7), and proceed to show that it is compatible
with the defining relations.

\medskip

{\bf 2.2. Multiplication table revisited.} We will first formally define
the action of $l_{\sigma}$ upon generators $\mu (\tau )$
sending them to some elements of $\Cal{M}_S$.
If $\sigma$ does not
break $\tau$, we put
$$
l_{\sigma}\mu (\tau )= 0.
\eqno(2.3)
$$
\smallskip

If $\sigma$ breaks $\tau$ at a vertex $v$, we put
$$
l_{\sigma}\mu (\tau )= \mu (\sigma*\tau ).
\eqno(2.4)
$$
\smallskip

Finally, if $\sigma$ breaks $\tau$ at an edge $e$, we
choose an arbitrary allowed quadruple of flags $I,J,K,L$
as in 1.3.1, Case (i), and put in the notations explained above
$$
l_{\sigma}\mu (\tau ) =
-\sum_{\alpha :\,IJ\alpha e_1} \mu (\tau (\alpha ))-
\sum_{\beta :\,e_2\beta KL} \mu (\tau (\beta))\, .
\eqno(2.5)
$$
\medskip

\proclaim{\quad 2.3. Theorem}  Formulas (2.3)--(2.5)
induce a well defined action of  $H^*_S$ upon $H_{*S}$.
\endproclaim

\medskip

\proclaim{\quad 2.3.1. Corollary} (a) This action
makes of $H_{*S}$ a free module of rank one over
$H^*_S.$

\smallskip

(b) The standard relations (1.10) span all linear relations
between the classes of good monomials $m(\tau )$ in $H^*_S$.
\endproclaim

\smallskip

{\bf Deduction of the Corollary.} Since the relations
between the classes of $\mu (\tau )$ correspond to a part of linear relations between good monomials, we have a surjective $k$--linear map
$$
s:\, H_{*S}\to H^*_S:\ [\mu (\tau )]\mapsto [m(\tau )].
$$
Here we set 
$[m (\tau )]:= m (\tau )\,\roman{mod}\,I_S$.

\smallskip

On the other hand, denoting by $\bold{1}\in H_{*S}$ the class
of one--vertex $S$--tree, thanks to Theorem 2.3 
we have another linear map
$$
t:\, H^*_S\to H{*_S}:\ [m (\tau )]\mapsto [m (\tau )]\,\bold{1}.
$$
Induction on the number of edges of $\tau$ using (2.4)
shows that $[m (\tau )]\,\bold{1}=[\mu (\tau )].$
Therefore $s$ and $t$ are mutually inverse, so that
$\bold{1}$ is a free generator of $H{*_S}$ as  
 an $H^*_S$--module. This completes the deduction.

\smallskip

We now turn to the proof of the theorem which is given
in the subsections 2.3.2 -- 2.3.6. It consists of a 
chain of sometimes tedious checks, and the reader may prefer
to skip it.

\medskip

{\bf 2.3.2. Prescription (2.5) is well defined.} The right hand side of
(2.5) formally depends on the choice of an allowed quadruple
$I,J,K,L$. One can pass from one allowed quadruple to another one
replacing each time only one flag so that all intermediate
quadruples will also be allowed (at each vertex, start with replacing
a flag carrying a white label by another flag also carrying a white label).

\smallskip

We will now take the difference
of the right hand sides of (2.5) written 
 for quadruples $I,J,K,L$ and $I',J,K,L.$ The second sums will cancel.
In the first sums, the terms corresponding to $\alpha$'s
with $II'J\alpha e_1$ will also cancel. The remaining terms
are
$$
-\sum_{\alpha:\,IJ\alpha I'e_1} \mu (\tau (\alpha ))
+ \sum_{\alpha:\,I'J\alpha Ie_1} \mu (\tau (\alpha ))\ \roman{mod}\, J_S.
\eqno(2.6)
$$
Clearly, (2.6) is an element of the form (2.1), namely
$r(\tau, v_1;I',J,I,e_1).$

\smallskip

One can similarly treat a replacement of $J$ by $J'$.

\medskip

{\bf 2.3.3. Multiplication by $l_{\sigma}$ maps $J_S$ into itself.}
As in (2.1), choose a relation  $r(\tau ,v;I,J,K,L)$
and a 2--partition $\sigma$.  We want to show that
$l_{\sigma} r(\tau ,v;I,J,K,L)$ belongs to $J_S$.
Consider in turn the following subcases.

\medskip

\centerline{\it (i) $\sigma$ does not break $\tau$}

\smallskip

If $\sigma$ does not break $\tau$, it does not break
any of $\tau (\alpha )$, hence 
$l_{\sigma} r(\tau ,v;I,J,K,L)=0$ in view of (2.3).

\medskip

\centerline{\it (ii) $\sigma$ breaks $\tau$ at a vertex different from $v$}

\smallskip

In this case we obviously have 
$$
l_{\sigma} r(\tau ,v;I,J,K,L)= r(\sigma*\tau ,v;I,J,K,L).
$$

\medskip

\centerline{\it (iii) $\sigma$ breaks $\tau$ at an edge not incident to $v$}

\smallskip

In this case we
calculate $l_{\sigma}\tau (\alpha )$
using one of the formulas (2.5) which involve
surgery only at the ends of $e$ so that it does
not interact with the partitions $\alpha$ in (2.1).
Hence we get a sum of the elements of the type (2.1)
written for various $\tau (\alpha ).$

\smallskip

The remaining two cases are more difficult.

\medskip

\centerline{\it (iv) $\sigma$ breaks $\tau$ at an edge $e$ incident to $v$}

\smallskip

We will start with some notation. Let $v_0$ be the vertex of
$e$ distinct from $v$. Denote by $e_0$, resp. $e_1$,
the flags of $e$ at $v_0$, resp. $v$. For each $\alpha$
appearing in (2.1), the tree $\tau (\alpha )$ has an extra edge
which we denote $f=f_{\alpha}$; other edges, in particular, $e$,
 ``are'' the edges of $\tau$. Clearly, $e$ and $f_{\alpha}$ have in 
$\tau (\alpha )$
a common vertex $v$. 

\smallskip

In order to calculate $l_{\sigma} \mu (\tau (\alpha ))$ we
have to use a formula of the type (2.5) written for
the edge  $e$ of $\tau (\alpha ).$ Notation in (2.5)
conflicts with that in (2.1), so we will rewrite (2.5) as:
$$
l_{\sigma}\mu (\tau (\alpha )) =
-\sum_{\beta :\,I_{\alpha}J_{\alpha}\beta e_0} \mu (\tau (\alpha )(\beta))-
\sum_{\gamma :\,e_1\gamma K_{\alpha}L_{\alpha}} \mu (\tau (\alpha)(\gamma))\, .
\eqno(2.7)
$$
Here $I_{\alpha},J_{\alpha},K_{\alpha},L_{\alpha}$ can
be chosen
depending on $\alpha$, thanks to 2.3.2. Recall from 1.3.1
that these flags must form an allowed quadruple, and 
$$
I_{\alpha},J_{\alpha}\in F_{\tau (\alpha )}(v_0)=F_{\tau} (v_0),\quad
K_{\alpha},L_{\alpha}\in F_{\tau (\alpha )}(v).
\eqno(2.8)
$$
Moreover, $\beta$, resp. $\gamma$, run over 2--partitions of
$F_{\tau (\alpha )}(v_0)=F_{\tau} (v_0)$, resp. of $F_{\tau (\alpha )}(v)$,
whereas $\alpha$ runs over partitions of $F_{\tau}(v).$

\smallskip

With this notation, we have
$$
l_{\sigma} r(\tau,v;I,J,K,L)=
 \sum_{\alpha} \epsilon (\alpha ;I,J,K,L)\, l_{\sigma}\mu (\tau (\alpha )) =
$$
$$
-\sum_{\alpha} \epsilon (\alpha ;I,J,K,L)\,\left(
\sum_{\beta :\,I_{\alpha}J_{\alpha}\beta e_0} \mu (\tau (\alpha )(\beta))+
\sum_{\gamma :\,e_1\gamma K_{\alpha}L_{\alpha}} \mu (\tau (\alpha)(\gamma))
\right) .
\eqno(2.9)
$$
Now, in (2.9) $\alpha$ and $\beta$ occur at different vertices of $\tau$.
Hence $\mu (\tau (\alpha )(\beta ))=\mu (\tau (\beta )(\alpha ))$,
and we can interchange summation on $\alpha$ and $\beta$
if we choose $I_{\alpha},J_{\alpha}$ independent of $\alpha$.
This is possible: namely, choose as $J_{\alpha}=J_0$
an arbitrary flag in $F_{\tau}(v_0)\setminus \{e_1\}$
whose branch carries a white label, and for $I_{\alpha}=I_0$
choose any other flag in this set. The condition
on $J_0$ will make a quadruple $I_0,J_0,K_{\alpha},L_{\alpha}$
allowed if  $L_{\alpha}$ carries a white label as well.
We will care about it later.

\smallskip

Meanwhile we can rewrite the first half of (2.9):
$$
-\sum_{\beta:\,I_0J_0\beta e_0} \sum_{\alpha}
\epsilon (\alpha ;I,J,K,L)\,
\mu (\tau (\beta )(\alpha))=
-\sum_{\beta:\,I_0J_0\beta e_0}
r(\tau (\beta ),v;I,J,K,L)\, \in J_S.
\eqno(2.10)
$$
Hence it remains to show that, after an appropriate choice of
$K_{\alpha},L_{\alpha}$, 
$$
(?)\quad \sum_{\alpha} \epsilon (\alpha ;I,J,K,L)\,
\sum_{\gamma :\,e_1\gamma K_{\alpha}L_{\alpha}} \mu (\tau (\alpha)(\gamma))
\in J_S.
\eqno(2.11)
$$

{\it Choice of $K_{\alpha},L_{\alpha}$.} For $L_{\alpha}$
we will always choose the flag of the new edge $f_{\alpha}$
incident to $v$. This choice assures that $I_0,J_0,K_{\alpha},
L_{\alpha}$ will form an allowed quadruple.

\smallskip
The choice of $K_{\alpha}$ will additionally depend
on the mutual position of $e$ and of flags $I,J,K,L$.
Notice that in our calculation this set of data is fixed: it 
determines the left hand side of (2.9).

\smallskip

There are five logical possibilities: one or none of the flags
$I,J,K,L$ coincides with $e_v$, the flag of $e$
incident to $v$. By symmetry, it suffices to consider
three options: none is $e_v$, $I=e_v$, or $J=e_v.$ 
We need only look at those $\alpha$ for which
$IJ\alpha KL$ or $KJ\alpha IL$.
\smallskip

Here is the list of choices.

\smallskip

If none of $I,J,K,L$ is $e_v$, we put $K_{\alpha}=I$
for $\alpha Ie_v$ and $K_{\alpha}=K$
for $\alpha Ke_v.$

\smallskip

If $I=e_v$, we put $K_{\alpha}=J$
for $\alpha Je_v$ and $K_{\alpha}=L$
for $\alpha Le_v.$

\smallskip

Finally, if $J=e_v$, we again put $K_{\alpha}=I$
for $\alpha Ie_v$ and $K_{\alpha}=K$
for $\alpha Ke_v.$

\medskip

Now the left hand side of (2.11) is well defined, and we
will start rewriting it as a sum of elements
of $J_S$. 

\smallskip

First of all, each term $\mu (\tau (\alpha) (\gamma ))$
can be uniquely rewritten as 
$\mu (\tau (\gamma') (\alpha' ))$ creating two new edges
of $\tau (\alpha) (\gamma )$ in reverse order.

\smallskip

More formally,  let $\tau (\gamma')$ be the result of
contracting the edge $f=f_{\alpha}$ in $\tau (\alpha) (\gamma )$.
It can be obtained from $\tau$ by inserting a new edge $g$
breaking the set of branches at the vertex $v$ into two parts. 
This edge has a common vertex with $e$.
Denote by $w=w_{\gamma'}$ another vertex of $g$.
The tree $\tau (\gamma') (\alpha' )$ 
is obtained from $\tau (\gamma')$ by inserting an edge
breaking the branches at $w$ into two parts.

\smallskip

Second, we will show that the sign 
$\epsilon (\alpha ;I,J,K,L)$ in (2.11)
can be rewritten as $\epsilon (\alpha' ;I',J',K',L')$
where the flags $I',J',K',L'$ at $w$ will depend
on $\gamma^{\prime}$ but not on $\alpha'$ and form
an allowed quadruple.

\smallskip

Then (2.11) can be rewritten as
$$
\sum_{\gamma'}\sum_{\alpha'}\epsilon (\alpha' ;I',J',K',L')\,
\mu (\tau (\gamma') (\alpha' )).
\eqno(2.12)
$$
It will remain only to convince ourselves that
$\alpha'$ runs over all 2--partitions of $F_{\tau (\gamma')}(w_{\gamma'}).$

\smallskip

We will dogmatically describe the choices involved.
We checked that everything fits together by drawing twelve
diagrams exhausting all possible mutual positions
of $I,J,K,L,e,f,g$ in the trees 
$\tau (\alpha) (\gamma ) = \tau (\gamma') (\alpha' ).$

\smallskip

The diagrams show that $\gamma'$ breaks $I,J,K,L$
in such a way that either all four flags get in the same part,
or one of them is taken apart. In the first case
we simply choose $I'=I,J'=j,K'=K,L'=L.$ In the second
case we replace the flag that is taken apart
by the appropriate flag belonging to the edge $g$.
One easily sees that the resulting quadruples are allowed.

\medskip

\centerline{\it (v) $\sigma$ breaks $\tau$ at $v$}

\smallskip

Again, we will consider several cases depending on how
$\sigma$ breaks $I,J,K,L$. There are four basic
options: $\sigma IJKL$, $I\sigma JKL$, $J\sigma IKL$,
$IJ\sigma KL$. All other possibilities either
can be obtained from one of these by exchanging
$J$ and $L$, $I$ and $K$, or else 
refer to those $\sigma$ for which
$l_{\sigma} r(\tau ,v;I,J,K,L)=0.$ The latter happens when
$IK\sigma JL$, because then $\sigma$ cannot break any
$\tau (\alpha )$ with either $IJ\alpha KL$ or $IJ\alpha KL$.

\smallskip

Assume first that one of the options $\sigma IJKL$, $I\sigma JKL$, 
$J\sigma IKL$
holds. Denote by $w$ the vertex of $\sigma *\tau$ incident
to three or more of the flags $I,J,K,L$.

\smallskip

As in (2.12), we can rewrite the nonvanishing terms of 
$l_{\sigma} r(\tau ,v;I,J,K,L)$
as follows:
$$
\epsilon (\alpha ;I,J,K,L)\, l_{\sigma}\mu (\tau (\alpha ))=
\epsilon (\alpha' ;I',J',K',L')\,\mu ((\sigma *\tau)(\alpha')).
\eqno(2.13)
$$
Here $\alpha'$ is a 2--partition of $F_{\sigma*\tau}(w)$,
and $I',J',K',L'$ are defined by the following prescription.
If $I,J,K,L$ are all contained in the same part of $\alpha$,
$I'=I,J'=J,K'=K,L'=L.$ If one of the flags is taken apart,
it must be replaced by the flag of $e$ incident to $w$.

\smallskip

A straightforward check then shows that  the summation over
all $\alpha$ produces the same result as the summation over all $\alpha'$,
so that finally
$$
l_{\sigma} r(\tau,v;I,J,K,L)=r(\sigma*\tau,w;I',J',K',L').
\eqno(2.14)
$$  

\smallskip

It remains to consider the case $IJ\sigma KL.$
In this case, among the summands of $r(\tau ,v;I,J,K,L)$
exactly one has the property $\tau (\alpha_0)=\sigma*\tau$.
To multiply it by $l_{\sigma}$, we use the formula
(2.5), in which we have to replace $\tau$ by $\tau (\alpha_0)$
and rename the summation variables, say, to $\beta ,\gamma$.
However, we  will use our current $I,J,K,L$ in the same sense as in (2.5),
so that finally we get
$$
l_{\sigma}\mu (\tau (\alpha_0))=
-\sum_{\beta :\,IJ\beta e_1} \mu ((\sigma*\tau) (\beta ))-
\sum_{\gamma :\,e_2\gamma KL} \mu ((\sigma*\tau )(\gamma ))\, .
\eqno(2.15)
$$
On the other hand, the nonvanishing terms $l_{\sigma}\mu (\tau (\alpha ))$ with $\alpha\ne\alpha_0$
can be represented as $\mu ((\sigma*\tau ) (\delta ))$ where $\delta$
breaks $\sigma*\tau $ at one of the vertices of $e$.
A contemplation shows that these terms exactly cancel (2.15).

\medskip

{\bf 2.3.4. Multiplications by $l_{\sigma}$ pairwise commute $\roman{mod}\,J_S.$} We want to prove that
$$
l_{\sigma_1}(l_{\sigma_2}\mu (\tau ))\equiv
l_{\sigma_2}(l_{\sigma_1}\mu (\tau ))\,\roman{mod}\,J_S.
\eqno(2.16)
$$
There is a long list of subcases that have to be treated
separately: ${\sigma_1}$, resp. $\sigma_2$, can break $\tau$
at a vertex $v_1$, resp. $v_2$, an edge $e_1$, resp. $e_2$,
or not to break $\tau$. For each of the possible 
breaking combinations, if the surgery loci are not incident,
say, $v_1$ is not a vertex of $e_2$, the reasoning
is straightforward, but it
becomes more tedious otherwise.

\smallskip

In the section A.3 of the Appendix to [LoMa], we displayed
the relevant calculations for the case when
the painted set $S$ has exactly two white labels.
In this case edges of $\tau$ have to form a linear chain,
with two white labels attached at the respective end vertices.

\smallskip

The total number of white labels influences everything in our
definitions: the total supply of painted stable trees,
the formulas for multiplication, and the list of the generators
of $I_S$ and $J_S$. Therefore generally we cannot simply refer
to the Appendix in [LoMa]. However, a check of pairwise commutativity 
of $l_{\sigma}$'s does
allow such a reduction.

\smallskip

Let us illustrate this in the case numbered $(ii)(ii')$ in A.3, [LoMa]: 
${\sigma_1}$, resp. $\sigma_2$, breaks $\tau$ at
the edges $e_1$, resp. $e_2$, having a common vertex.
Let us replace $\tau$ by the linear tree $\tau_0$ with  edges
$e_1,e_2$. Besides halves of these edges, we endow $\tau_0$
with tails which are branches of $\tau$ incident to
the respective three vertices. At the two end vertices of $\tau_0$,
choose two branches
of $\tau$  which carry white labels
and declare the respective tails of $\tau_0$ white.
Declare other tails black. 
Replace the 2--partitions $\sigma_1, \sigma_2$, by the induced
partitions of tails of $\tau_0$. 

\smallskip

Now, the calculations in A.3 referring to $\tau_0$
are compatible with the respective calculations for
$\tau$, because the formula (2.5) is ``local''.
They allow less freedom with respect to the
choice of $I,J,K,L$, but
this is inessential, because from 2.3.2 and 2.3.3
it follows that any choice will do. 

\smallskip

Finally,
with the choices made in [LoMa], we get (2.16) as an exact
equality.

\medskip

We have now proved that $\Cal{R}_S$ acts upon $H_{*S}.$
It remains to check that $I_S$ annihilates $H_{*S}.$

\medskip

{\bf 2.3.5. Elements $R_{ijkl}$ annihilate $H_{*S}$.}
The case when all labels are white was treated in
[Ma1], Chapter III, 4.8.4. With very minor
additional precautions, the treatment can be repeated
for general painted $S$.

\smallskip

The reason is that, although we have a short supply of
painted stable 2--partitions $\sigma$'s if not all labels are white,
nevertheless for an allowed quadruple $i,j,k,l\subset S$,
the relation $R_{ijkl}$ involves with non--zero coefficients 
exactly those $\sigma$'s which are simply stable.
Painted stable trees also are simply stable.
Finally, the multiplication formulas (2.5) written
for allowed quadruples $I,J,K,L$ constitute a part
of the multiplication formulas in the unpainted case.

\smallskip

In [Ma1], III.4.8.4, we show that $R_{ijkl}[\mu (\tau )]=0,$
again for different reasons in different combinatorial
situations. When no choices are involved in a
calculation, it is valid in the painted case as well.
If a choice of $I,J,K,L$ is involved,
one can check that for an allowed $i,j,k,l$,
the choices made in [Ma1] are allowed as well.

\medskip

{\bf 2.3.6. Elements $R_{\sigma_1\sigma_2}$ annihilate $H_{*S}$.}
These elements are $l_{\sigma_1}l_{\sigma_2}$ such
that $\sigma_1$ and $\sigma_2$ do not break each other.
If one of them, say, $\sigma_1$, does not break $\tau$ either,
then because of commutativity $l_{\sigma_1}l_{\sigma_2}\mu (\tau )=0.$

\smallskip

If both of them break $\tau$, they have to break $\tau$ at one
and the same vertex $v$. It is then easy to check that
$l_{\sigma_1}l_{\sigma_2}\mu (\tau )=0.$

\newpage

\centerline{\bf \S 3. Cohomology and homology as functors on trees}

\medskip

{\bf 3.1. Groups $H^*(\tau )$ and $H_*(\tau )$.} Let $\tau$ be
a stably painted tree. We extend the painting of tails of $\tau$
to the painting of all flags declaring all halves of the edges
white. Finally, we put
$$
H^*(\tau ):=\otimes_{v\in V_{\tau}} H^*_{F_{\tau}(v)},\quad
H_*(\tau ):=\otimes_{v\in V_{\tau}} H_{*F_{\tau}(v)} \,.
\eqno(3.1)
$$
For one--vertex trees, we recover  the rings/modules (1.3), (2.2) 
with which we
worked in the earlier sections. The tensor products (over $k$)
of families of objects labeled by finite sets, as in (3.1),
are defined generally in symmetric monoidal categories
and allow one to make explicit the functorial properties with respect
to the maps of the index sets. In particular, $H^*(\tau )$, resp. $H_*(\tau )$,
has a natural structure of $k$--algebra, resp. $H^*(\tau )$--module, with 
compositions defined componentwise.

\smallskip

The main result of this section is:

\medskip

\proclaim{\quad 3.2. Theorem} For any contraction morphism
$f:\, \tau'\to\tau$ of stably painted trees, 
one can define a ring homomorphism
$$
f^*:\,H^*(\tau )\to H^*(\tau' )
\eqno (3.2)
$$
and a compatible homomorphism of modules
$$
f_*:\,H_*(\tau' )\to H_*(\tau )
\eqno (3.3)
$$
which makes $H^*$ and $H_*$ functors on the category 
of trees and contractions.

\smallskip

Functors $H^*,H_*$ are well determined by their actions
on isomorphisms (obvious) and
on the morphisms contracting a unique edge of an one--edge
tree. These restrictions are explicitly described below.
\endproclaim

\medskip

{\bf Proof.} Let $S$ be a painted stable set identified
with an one--vertex tree with flags $S$, $\sigma$
an one--edge painted stable $S$--tree identified with 
a 2--partition $S=S_1\cup S_2$. We denote by $e$ the
edge of $\sigma$, and by $v_1,v_2$ 
its vertices carrying flags $S_1,S_2,$ 
and by $e_1,e_2$ the respective halves of $e$. 
Denote by $f_e$ the contraction morphism. We want first of all to
define a ring homomorphism $f^*_e: H^*_S\to H^*_{S_1\cup\{e_1\}}
\otimes H^*_{S_2\cup\{e_2\}}$.

\smallskip

We start with describing its (depending on some auxiliary choices) lift
$\phi_e^*:\,\Cal{R}_S\to \Cal{R}_{S_1\cup\{e_1\}}\otimes \Cal{R}_{S_2\cup\{e_2\}}.$
Put $\phi_e^*(l_\rho )=0$ if $\rho$ does not break $\sigma$.
If $\rho$ breaks $\sigma$ at $v_1$, it defines a painted stable
2--partition $\rho_1$ of $S_1$ and the respective new edge
$f$ of $\rho*\sigma$. Let $S_1'$ denote the set of flags of this tree
at the common vertex of $e$ and $f$, and $S_1''$ the
flags at the other vertex of $f$. Denote by $\rho_1$
the one edge $(S_1\cup\{e_1\})$--tree with tails $S_1'\cup\{e_1\}$, $S_1''$.
Similarly, if  $\rho$ breaks $\sigma$ at $v_2$ define
an one--edge $(S_2\cup\{e_2\})$--tree $\rho_2.$ Put
$$
\phi_e^*(l_{\rho}):=l_{\rho_1}\otimes 1,\ \roman{resp.}\ 1\otimes l_{\rho_2}
\in \Cal{R}_{S_{1\cup\{e_1\}}}\otimes \Cal{R}_{S_{2\cup\{e_2\}}}.
\eqno(3.4)
$$
It remains to define $\phi_e^*$ in the case when $\rho =\sigma$,
or equivalently, when $\rho$ breaks $\sigma$ at $e$.
In this case we choose an allowed quadruple of flags
$i,j,k,l$ such that $i,j\in S_1,\,k,l\in S_2,$ and define, similarly to
(1.7), (2.5),
$$
\phi_e^*(l_{\sigma})= -\sum_{\alpha :\,ij\alpha e_1}l_{\alpha}\otimes 1
-\sum_{\beta :\,ij\beta e_1} 1\otimes l_{\beta} \,\in 
 \Cal{R}_{S_{1\cup\{e_1\}}}\otimes \Cal{R}_{S_{2\cup\{e_2\}}}.
\eqno(3.5)
$$ 
Notice right away that another choice of $i,j,k,l$ will 
not change the result modulo $I_{S_{1\cup\{e_1\}}}\otimes
\Cal{R}_{S_{2\cup\{e_2\}}}\,+ 
\Cal{R}_{S_{1\cup\{e_1\}}}\otimes
I_{S_{2\cup\{e_2\}}}.$ The same reasoning as in 2.3.2
shows this: we consider a replacement of one flag, say, $i'$, and then
the respective difference as in (2.6)
can be identified with
$r(\alpha ,v_1;i',j,i,e_1)\otimes 1.$

\medskip

\proclaim{\quad 3.2.1. Lemma} (a) There is a well defined ring homomorphism
$$
f^*_e: H^*_S\to H^*_{S_1\cup\{e_1\}} \otimes H^*_{S_2\cup\{e_2\}},\quad
[l_{\rho}]\mapsto [\phi_e^*(l_{\rho})]\,.
\eqno(3.6)
$$
(ii) $f_e^*$ is surjective. More precisely, 
let $\tau_1$, resp. $\tau_2$, be a painted stable
$(S_1\cup\{e_1\})$--tree, resp. $(S_2\cup\{e_2\})$--tree.
Denote by $\tau_1\bullet \tau_2$ the $S$--tree
obtained by gluing $e_1$ to $e_2$. Then
$$
f_e^*([m(\tau_1\bullet \tau_2)])=[m(\tau_1 )]\otimes [m(\tau_2 )].
\eqno(3.7)
$$
\endproclaim

\medskip

{\bf Proof of Lemma 3.2.1.}  We must show that
$\phi^*_e$ maps the generators (1.1) and (1.2) of $I_S$
into $I^*_{S_1\cup\{e_1\}} \otimes \Cal{R}^*_{S_2\cup\{e_2\}}\,+
\Cal{R}^*_{S_1\cup\{e_1\}} \otimes I^*_{S_2\cup\{e_2\}}.$

\smallskip

Start with a generator $R_{ijkl}$ which we will now endow
with a superscript $S$. Our reasoning depends on the mutual
position of $i,j,k,l$ and $\sigma.$

\smallskip

{\it Case $ijkl\sigma$.} If, say, $i,j,k,l\in S_2$, we have
$$
\phi_e^*(R_{ijkl}^S)=1\otimes R^{S_2\cup\{e_2\}}_{ijkl}\,.
$$
\smallskip

{\it Case $j\sigma ikl$.} If $i,k,l\in S_2$, we have
$$
\phi_e^*(R_{ijkl}^S)=1\otimes R^{S_2\cup\{e_2\}}_{ie_2kl}\,.
$$
\smallskip

{\it Case $i\sigma jkl$.} If $j,k,l\in S_2$, we have
$$
\phi_e^*(R_{ijkl}^S)=1\otimes R^{S_2\cup\{e_2\}}_{e_2jkl}\,.
$$
\smallskip

{\it Case $ik\sigma jl$.} In this case 
$$
\phi_e^*(R_{ijkl}^S)=0\,.
$$
\smallskip

{\it Case $ij\sigma kl$.} In this case one summand
of $R_{ijkl}^S$ is $l_{\sigma}.$ Applying $\phi^*_e$
to it with the same choice of $i,j,k,l$  we get the expression (3.5).
Applying $\phi^*_e$ to all other terms of
$R_{ijkl}^S$ and using (3.4), we get the same terms as in
(3.5), with reverse signs. Hence finally
$$
\phi_e^*(R_{ijkl}^S)=0\,.
$$

\smallskip

Up to obvious symmetries, we have exhausted all possible
alternatives. It remains to show that
$$
\phi_e^*(l_{\sigma_1})\phi_e^*(l_{\sigma_2})\in 
I_{S_1\cup\{e_1\}} \otimes \Cal{R}_{S_2\cup\{e_2\}}\,+
\Cal{R}_{S_1\cup\{e_1\}} \otimes I_{S_2\cup\{e_2\}}
$$
if $\delta (\sigma_1,\sigma_2)=2.$ We leave this 
as an exercise to the reader.

\smallskip

We have thus  established that $\phi^*_e$ induces
a ring homomorphism $f_e^*$. Representing
$m(\tau_1\bullet \tau_2)$ as the product
of generators corresponding to the edges, and applying (3.4),
we get (3.7). This shows the surjectivity
of $f^*_e$ and completes the proof of the lemma.

\medskip

We now continue the proof of the Theorem 3.2 and turn to
$H_*$. Motivated by (3.7) (and, of course, by algebraic geometry),
we define (in the notations of Lemma 3.2.1 (ii)):
$$
f_{e*}:\,H_{*S_1\cup\{e_1\}} \otimes H_{*S_2\cup\{e_2\}}
\to H_{*S},\quad [\mu (\tau_1)]\otimes [\mu (\tau_2)]
\mapsto [\mu (\tau_1\bullet\tau_2)].
\eqno(3.8)
$$
Compatibility of this prescription with the defining relations
of $J_S$ (cf. (2.1)) is straightforward.

\smallskip

Since $[\mu (\tau )]=[m(\tau )]\,\bold{1}$ (cf. the proof of 2.3.1),
combining (3.7) and (3.8) we obtain
$$
f_{e*}(f_e^*([m(\tau_1\bullet\tau_2)])\,(\bold{1}\otimes\bold{1}))=
[m(\tau_1\bullet\tau_2)]\,\bold{1}.
$$
This is the key special case of the compatibility of
$f_e^*$ and $f_{e*}$ (``projection formula'' of algebraic geometry).
The general case (obtained by replacing $\bold{1}\otimes\bold{1}$
with any element of $H_*\otimes H_*$) follows from it formally,
because $\bold{1}$ is a free generator of $H_*$ over the respective $H^*$.

\smallskip

We can now define (3.2) and (3.3) by decomposing
$f$ into a product of edge contractions and a final isomorphism
and using (3.7), (3.8). Edge contractions commute
in an intuitively evident sense, so that checking that the result
is independent of the decomposition involves only
a careful bookkeeping. We leave this as an exercise
to the reader.

\bigskip
\centerline{\bf \S 4. $\Cal{L}$--algebras}

\medskip

In this section we will define several versions of
the notion of ``algebra over extended modular operad of genus zero''.
The basic definition we start with is modeled upon the 
functorial treatment introduced
in [KoMa], sec. 6. It involves unoriented painted stable trees and
(combinatorial) cohomology algebras, and it
produces a version of {\it cyclic
algebras} in the sense of Getzler and Kapranov. Tensor product
of algebras appears in the most straightforward way in this construction.

\smallskip

Other versions involve oriented trees and/or homology algebras.
The generic name we use for  all these versions is $\Cal{L}$--{\it algebras.}

\smallskip

Finally, we will explain how the classical operads appear in this context.

\medskip

{\bf 4.1. From trees to tensors: unoriented case.}
Let $\Cal{A}=\{T;F, (\,,)\}$ be a triple consisting
of two free (or projective) finite rank $\bold{Z}_2$--graded $k$--modules
$T, F$ and a (super)symmetric scalar product
$(\,,)$ on $F$. We assume that $(\,,)$ induces an isomorphism
$F\to F^t$ where $F^t$ is the graded dual to $F.$
We denote by $\Delta\in F\otimes F$ the respective Casimir element.

\smallskip

Let $\tau$ be a painted stable tree. Denote by $W_{\tau}$, resp. $B_{\tau}$,
 the set of
its white, resp. black flags (recall that halves of edges are all white.)
If $v$ is a vertex of $\tau$, $W_{\tau}(v)$, resp. $B_{\tau}(v)$,
denotes the set of white, resp. black, flags incident to this vertex.
Put 
$$
\Cal{A}(\tau):= T^{\otimes B_{\tau}}\otimes
F^{\otimes W_{\tau}}=
\otimes_{v\in V_{\tau}} T^{\otimes B_{\tau}(v)}\otimes
F^{\otimes W_{\tau}(v)}\,.
\eqno(4.1)
$$
Thus, we attach the ``white space'' $F$ at each white flag,
the ``black space'' $T$ at each black flag, and take
the tensor product of all spaces.

\smallskip

Let now $f:\,\tau\to\sigma$ be a contraction. It identifies
$F_{\sigma}$ with a subset $F^s_{\tau}$ of non--contracted
flags. The complement consist of the halves
of contracted edges $E^c_{\tau}.$ We can define a natural map
$$
f_{\Cal{A}}^*:\,\Cal{A}(\sigma )\to \Cal{A}(\tau ) 
\eqno(4.2)
$$
which tautologically identifies spaces attached to the
flags which are paired by $f$ and then tensor multiplies
the result by 
$\Delta^{\otimes E^c_{\tau}}\in 
(F\otimes F)^{\otimes E^c_{\tau}}$.

\medskip

\proclaim{\quad 4.1.1. Proposition} 
In this way $\tau\mapsto \Cal{A} (\tau )$ becomes a contravariant functor
on the category of stably painted trees and contractions.
\endproclaim

\smallskip

This is straightforward.

\medskip

\proclaim{\quad 4.2. Definition} The structure of a cyclic $\Cal{L}$--algebra
upon $\Cal{A}=\{T;F, (,)\}$ is a morphism
of functors compatible with gluing
$$
I:\ \Cal{A}\to H^*.
\eqno(4.3)
$$
In other words, it consists of a family of maps indexed by stably painted trees
$$
I(\tau ):\ T^{\otimes B_{\tau}}\otimes
F^{\otimes W_{\tau}}\to H^*(\tau )
\eqno(4.4)
$$
such that for any contraction morphism $f:\,\tau\to\sigma $
we have
$$
I(\tau )\circ f^*_{\Cal{A}}=f^*\circ I(\sigma ):\
T^{\otimes B_{\sigma}}\otimes
F^{\otimes W_{\sigma}}\to H^*(\tau )
\eqno(4.5)
$$
and moreover
$$
I(\tau_1\bullet\tau_2)=I(\tau_1)\otimes I(\tau_2):\,
\Cal{A}(\tau_1)\otimes \Cal{A}(\tau_2)\to H^*(\tau_1)\otimes
H^*(\tau_2)=H^*(\tau_1\bullet\tau_2).
\eqno(4.6)
$$
\endproclaim

\smallskip

Restricting this definition to the case of pure white painted trees, we get
the structure on $F$ which was called
an operadic tree level
Cohomological Field Theory in [KoMa], Def. 6.10. This is the same as an algebra
over the cyclic operad $\{H_*(\overline{M}_{0,n+1},k)\}$.
The case of trees with exactly
two white flags (and orientation, cf. below) 
produces a structure which is essentially a representation
of a certain algebra $H_*T$ on $F$ studied
in [LoMa], sec. 3.3--3.6. Thus, the notion of an $\Cal{L}$--algebra combines
both these structures.

\medskip

{\bf 4.2.1. Tensor product of cyclic $\Cal{L}$--algebras.} Let
$\Cal{A}=\{T_1;F_1, (\,,)_1\},\, I_1$, resp. $\Cal{B}=\{T_2;F_2, (\,,)_2\},\,I_2$, be two cyclic $\Cal{L}$--algebras. Put
$$
\Cal{A}\otimes\Cal{B}:=\{T_1\otimes T_2;\,
F_1 \otimes F_2, (\,,)_1 \otimes (\,,)_2\},
$$
$$
I_{\Cal{A}\otimes\Cal{B}}(\tau ):=M\circ(I_1\otimes I_2)(\tau ):\,
$$
where $M:\,H^*(\tau )\otimes H^*(\tau )\to H^*(\tau )$
is the ring multiplication. Compatibility with
(4.5), (4.6) is straightforward.

\medskip

{\bf 4.2.2. Economy class descriptions of cyclic algebras.}
Let us consider only maps (4.4) indexed by
one--vertex trees, with  a painted set of flags $S=W\,\cup\,B$:
$$
I_S :\, T^{\otimes B}\,\otimes\, F^{\otimes W}\to H^*_S .
\eqno(4.7)
$$
The axiom (4.7) puts certain restrictions upon this family
of maps. Namely, for any painted stable partition $S=S_1\,\cup\,S_2$,
$S_i=W_i\,\cup\,B_i$, the following diagram must be commutative:
$$
\CD
T^{\otimes B}\,\otimes\, F^{\otimes W} @>I_S>> H^*_S \\
@VS\circ (\otimes \Delta )VV          @VVf^*V \\
T^{\otimes B_1}\,\otimes\, F^{\otimes W_1\cup\{e_1\}}\otimes
T^{\otimes B_2}\,\otimes\, F^{\otimes W_2\cup\{e_2\}}
@>I_{S_1\cup\{e_1\}}\otimes I_{S_2\cup\{e_2\}}>> 
H^*_{S_1\cup\{e_1\}}\otimes H^*_{S_2\cup\{e_2\}} \\
\endCD
\eqno(4.8)
$$
Here $f^*$ is the map defined by (3.6), (3.7), and $S$ is
an obvious reshuffling of factors.

\smallskip

Moreover, any painted isomorphism $S\to S^{\prime}$
produces isomorphisms of the respective corners
of the diagram (4.8), and we require 
the resulting cubic diagram to be commutative.

\smallskip

One easily sees that, conversely, given $I_S$ satisfying
(4.8) and compatible with isomorphisms, we can uniquely
reconstruct the structure of cyclic $\Cal{L}$--algebra
on $\Cal{A}=\{T;F, (,)\}$.

\smallskip

We can further restrict ourselves to the subcategory of
painted sets $S_{n,m}$ in which $W$ (resp. $B$) consists of an initial
segment of the natural numbers
$\{1,\dots ,n\}$ (resp. an initial
segment of another copy of the natural numbers
$\{\bar{1},\dots ,\bar{m}\}$. Isomorphisms will reduce
to the actions of two copies of symmetric groups
$\bold{S}_n\times \bold{S}_m.$ The axioms
we get in this way are extensions of ones that are stated
in [KoMa], Definition 6.1.

\medskip

{\bf 4.3. From trees to tensors: oriented case.} An {\it orientation}
of a painted stable root tree $\tau$ is uniquely determined by a choice
of its {\it root}: an arbitrary white flag $r$.
In the geometric realization, we then orient $r$ away from its vertex,
and each other white flag in the direction of $r$. 
Put $W_{\tau}=W_{\tau}^+\cup W_{\tau}^-$
where $W_{\tau}^+$ consists of flags oriented towards their vertices
(``in''),
and $W_{\tau}^-$ of flags oriented in reverse direction (``out''.)
Notice that all black flags are oriented ``in''.

\smallskip

An oriented contraction of oriented trees is the same
as a usual contraction which sends root to root.

\smallskip

To produce tensors associated with oriented trees we start
with a pair of modules $\Cal{A}_0=\{T,\,F\}$ as in 4.1,
but do not assume a scalar product given.
We then denote by $\Delta\in F\otimes F^t$
the image of the identity morphism.

\smallskip

For an oriented tree $\tau$ we put
$$
\Cal{A}_0(\tau):= T^{\otimes B_{\tau}}\otimes
(F^t)^{\otimes W_{\tau}^-}\otimes F^{\otimes W_{\tau}^+}\,.
$$

Let now $f:\,\tau\to\sigma$ be an oriented contraction. 
As in 4.1, we can define a natural map
$$
f_{\Cal{A}_0}^*:\,\Cal{A}_0(\sigma )\to \Cal{A}_0(\tau ) 
\eqno(4.9)
$$
which tautologically identifies spaces attached to the
flags which are paired by $f$ and then tensor multiplies
the result by 
$\Delta^{\otimes E^c_{\tau}}\in 
(F \otimes F^t)^{\otimes E^c_{\tau}}$.

\medskip

\proclaim{\quad 4.3.1. Proposition} 
In this way $\tau\mapsto \Cal{A}_0 (\tau )$ becomes a contravariant functor
on the category of oriented stably painted trees and contractions.
\endproclaim

\medskip

\proclaim{\quad 4.4. Definition} The structure of an $\Cal{L}$--algebra
upon $\Cal{A}_0=\{T,\,F\}$ is a morphism
of functors compatible with gluing
$$
I:\ \Cal{A}_0\to H^*.
\eqno(4.10)
$$
\endproclaim

\medskip

{\bf 4.5. $\Cal{L}$--algebras in terms of combinatorial homology and correlators.}
We will now introduce one more version of $\Cal{L}$--algebras
involving homology and unoriented trees. It is this version which
is most convenient for the passage to the differential equations.

\smallskip

For a change, we will start with an economy class description.
We fix  $\Cal{A}_0=\{T,\,F\}$ as in 4.3. 

\medskip

{\bf 4.5.1. Data.} Consider one-vertex oriented stable painted trees.
Such a tree is the same as a finite set of its flags represented
in the form $S=(W\cup \{0\})\,\cup\,B$, $|W\,\cup\,B|\ge 2,$ where $0$ is the root,
$W\cup \{0\}$ (resp. $B$) is the subset of white (resp. black) flags. 
For each such $S$ 
we assume given an even  homomorphism of $k$--modules (correlator)
$$
C_S:\,H_{*S}\otimes T^{\otimes B}\otimes F^{\otimes W}\to F\, .
\eqno(4.11)
$$
This should be compared to (4.7): to pass from (4.7) to (4.11), we
make a partial dualization: cohomology becomes homology
and moves to the left hand side,
whereas the factor $F^t$ corresponding to the root
becomes $F$ and moves to the right hand side.

\smallskip

Once again, we can restrict ourselves by the sets
$A=\{0,1,\dots ,n\}, B=\{\overline{1},\dots ,\overline{m}\}$
taken from two disjoint copies of the set of nonnegative integers.

\smallskip

Since $H_{*S}$ is spanned by $[\mu (\tau )]$ where $\tau$ runs over
unoriented stable painted $S$--trees, to give (4.11) is the same
as to give a family of maps
$$
C_S([\mu (\tau )]):\,T^{\otimes B}\otimes F^{\otimes W}
\to F\, .
\eqno(4.12)
$$
Notice again that the orientation is encoded by $S$ which is given
together with the choice of (a label of the) root.
\smallskip

The intuitive meaning of (4.12) is this: given a rooted tree $\tau$,
we consider it as the flowchart of a computation.
Concretely, an input assigns  an element $f_a$
of $F$ (resp. $t_b$ of $T$) to each white non--root tail $a$ (resp. to each 
black tail $b$). Then the output is an element
$$
C_S([\mu (\tau )])(\otimes_b\,t_b\otimes_a\,f_a)\in F.
\eqno(4.13)
$$
We associate this output with the root of $\tau$.

\medskip

{\bf 4.5.2. Axioms.} (A) {\it Operadic relations.} If we graft two (or more) oriented trees  
by gluing root(s) to white tail(s), the resulting
tree must be the flowchart for the composition
of the respective functions.

\smallskip

Notice that since black tails are never glued, elements of
$T$ may be considered as ``parameters'', whereas the
actual arguments (and values) of the correlators belong
to $F$. 

\smallskip

(B) {\it Symmetry relations.} They say that the product 
$\bold{S}_W\times\bold{S}_B$
of the permutation groups of $W$ and $B$ acting upon
(4.11) in an evident way
leaves these maps invariant. Notice that
the action on $H_{*S}$ is also involved in this requirement.

\smallskip

(C) {\it Linear relations.} They assure us that all
linear relations between $[\mu (\tau )]$ in $H_{*S}$
imply the respective relations between correlators.
We can restrict ourselves by the basic relations (2.1),
or their cohomological version (1.10).
Namely, starting with an $S$--tree $\tau$, its vertex $v$,
and an allowed quadruple of flags $I,J,K,L$ at this vertex,
as in 1.4, we must have
$$
\sum_{\alpha} \epsilon (\alpha ;I,J,K,L)\,
C_S([\mu (\tau (\alpha ))]) =0
\eqno(4.14)
$$
where we sum over all 2--partitions of the set of flags at $v$,
and $\varepsilon = 1, -1,$ or 0, depending on
whether $IJ\alpha KL,\, KJ\alpha IL$ or none of it.
\medskip

{\bf 4.6. Top correlators.} Let $\tau$ be an one vertex tree
whose tails are labeled by $(W\,\cup\,\{0\})\,\cup\, B$. For
$W=\{1,\dots ,n+1\},\, B=\{\overline{1},\dots ,\overline{m}\}$, 
we will write
$c(t_1,\dots ,t_m,f_1,\dots ,f_{n+1})$ in place of (4.13)
and call these elements  {\it the top correlators}.
They are symmetric in $T$-- and $F$--arguments separately. 
We will formally extend them to polylinear functions
symmetric in all variables. All the operadic
relations follow by iteration from this symmetry and from
the following ones:
$$
c(t_1,\dots ,t_m,\,f_1, f_2,\dots ,f_n, c({t'}_1,\dots ,{t'}_q;
{f'}_1,\dots ,{f'}_{p+1}))=
$$
$$
c(t_1,\dots ,t_m,f_1,f_2,\dots ,f_n,{t'}_1,\dots ,{t'}_q,
{f'}_1,\dots ,{f'}_{p+1})
\eqno(4.15)
$$
Now fixing all the arguments except for
$f_{n+1}$, consider  $c(t_1,\dots ,t_m,f_1,\dots ,f_{n+1})$
as a function of $f_{n+1}$, that is, an element of
$\roman{End}\,F$ (even or odd, depending on the sum
of parities of the fixed arguments), and denote it by
$$
\langle t_1,\dots ,t_m,f_1,\dots ,f_{n}\rangle \,\in\, \roman{End}\,F\,.
$$

Again, we will allow arbitrary permutations of arguments
in this expression, followed by the standard sign change,
and call the resulting expressions {\it top matrix correlators}.
Then the grafting relations (4.15) turn into the
simple
{\it factorization relations} in $\roman{End}\,F$:
$$
\langle t_1,\dots ,t_m,f_1,f_2,\dots ,f_n,{t'}_1,\dots ,{t'}_q,
{f'}_1,\dots ,{f'}_{p}\rangle =
$$
$$
\langle t_1,\dots ,t_m,f_1,f_2,\dots ,f_n\rangle
\langle{t'}_1,\dots ,{t'}_q,
{f'}_1,\dots ,{f'}_{p}\rangle .
\eqno(4.16)
$$
\smallskip
We can now formulate the first result
describing $\Cal{L}$--algebras in terms of their
top matrix correlators. In order to make connection
with results of [LoMa], we introduce the following notation.
Choose bases of $F,\,T$ and  let their union,
which is a basis of $T\oplus F$, be $\{\Delta_i\}$,
for some index set $I.$ Fix a structure
of the (oriented, homological)
 $\Cal{L}$--algebra upon $(T,F)$ and consider its top matrix correlators
for the elements of the chosen basis.

\medskip

\proclaim{\quad 4.6.1. Theorem} (i) Top matrix correlators 
$$
\langle \Delta_{a_1}\dots\Delta_{a_n}\rangle \in \roman{End}\,F,\quad n\ge 2,
$$
are (super)symmetric in their arguments. 
Furthermore, for any allowed quadruple $i,j,k,l\in I$ 
$$
\sum_{\sigma:\,ij\sigma kl}\varepsilon (\sigma,(a_m))
\langle\prod_{m\in\sigma_1}\Delta_{a_m}\rangle\cdot
\langle\prod_{m\in\sigma_2}\Delta_{a_m}\rangle - 
\sum_{\sigma:\,kj\sigma il}\varepsilon (\sigma,(a_m))
\langle\prod_{m\in\sigma_1}\Delta_{a_m}\rangle\cdot
\langle\prod_{m\in\sigma_2}\Delta_{a_m}\rangle =0\,.
\eqno(4.17)
$$
Here $\sigma$ runs over ordered
2--partitions of $\{1,\dots ,n\}$.
We choose additionally an arbitrary ordering of both parts $\sigma_1,
\sigma_2$ determining the ordering of $\Delta$'s in the
angular brackets, and compensate this choice by the
$\pm 1$--factor $\varepsilon (\sigma,(a_k))$.

\medskip

(ii) Conversely, any family of elements 
$\langle\Delta_{a_1}\dots\Delta_{a_n}\rangle\in \roman{End}\,F$ defined for all $n\ge 2$ and
$(a_1,\dots ,a_n)\in I^n$ and satisfying the symmetry conditions and (4.17),
is the family of top matrix correlators of a unique
structure of $\Cal{L}$--algebra on $(T,F)$.
\endproclaim 

\smallskip
 
{\bf Proof.} (i) Both symmetry properties follow from the
definitions,
whereas (4.17) follows from (4.14), written
for one--vertex trees $\tau$ and transformed with the help
of (4.16).

\smallskip

(ii) If we know all top matrix correlators, we can recursively
determine maps (4.13) for all stable trees $\tau $ with
$e\ge 1$ edges. In fact, choose an end edge of $\tau$
and cut it, thus representing $\tau$ as the result
of grafting an one--vertex tree $\rho$ to a tree $\sigma$ with
$e-1$ edges. Choose the flag $a$ of $\sigma$ that gets
grafted as special input and, using partial
dualization, rewrite (4.13) as a map
$T^{\otimes B}\otimes F^{\otimes W}\to 
\roman{End}\,F$, that is, a matrix correlator
related to the tree $\sigma$. From the grafting
axiom (A) above it follows that (4.13) is uniquely
calculated from the latter matrix correlator
and the top matrix correlator for $\rho$,
via a generalization of the factorization relations
(4.16).

\smallskip

This shows that the $\Cal{L}$--structure on $(T,F)$ with given top
matrix correlators is unique, if it exists
at all. In order to prove existence, we must check
that all linear relations (4.14) follow from their
special cases (4.17). The argument is essentially the same
as in the proof on p. 462 of [LoMa]; the difference is
that in [LoMa] the relevant products of top correlators were linear
ordered, whereas here they are controlled by a tree, so a more careful
bookkeeping is required. We leave this as an exercise.

\bigskip

\centerline{\bf \S 5. Differential geometry and $\Cal{L}$--algebras}

\medskip

{\bf 5.1. Notation.}  Consider a structure of $\Cal{L}$--algebra
upon $(T,F)$ determined by its top matrix correlators
$\langle \Delta_{a_1}\dots\Delta_{a_n}\rangle $ as in Theorem 4.6.1 (i).
Denote by $(x^a)$ the linear coordinate system on $T\oplus F$
dual to $(\Delta_a).$ Put
$$
\Cal{B} =
\sum_{n=1}^{\infty}\sum_{(a_1,\dots ,a_n)}
\frac{x^{a_n}\dots x^{a_1}}{n!}\,
\langle \Delta_{a_1} \dots \Delta_{a_n}\rangle \in k[[x]]\otimes \roman{End}\,F.
\eqno(5.1)
$$
\smallskip

\proclaim{\quad 5.1.1. Theorem} a) We have
$$
d\Cal{B}\wedge d\Cal{B}=0.
\eqno(5.2)
$$

\smallskip

b) Conversely, let $\Cal{B}\in k[[x]]\otimes \roman{End}\,F$
be any even formal series without constant term satisfying
the equation (5.2). Define its symmetrized coefficients
$\Delta (a_1, \dots ,a_n),$ $a_i\in I,$ by the following
conditions: first,
$$
\Cal{B} =
\sum_{n=1}^{\infty}\sum_{(a_1,\dots ,a_n)}
\frac{x^{a_n}\dots x^{a_1}}{n!}\,\Delta (a_1,\dots ,a_n),
$$
second, the parity of $\Delta (a_1, \dots ,a_n)$ coincides with
that of $x^{a_n}\dots x^{a_1}$, and third,
$$
\Delta (a_{s(1)},\dots ,a_{s(n)})=\varepsilon (s,(a_i))\,\Delta (a_1,\dots ,a_n)\,
$$
for any permutation $s$ of $\{1,\dots ,n\}.$

\smallskip

Then there exists a unique structure of $\Cal{L}$--algebra
upon $(T,F)$ with top correlators 
$$
\langle \Delta_{a_1} \dots \Delta_{a_n}\rangle =
\Delta (a_1,\dots ,a_n).
$$
\endproclaim

\smallskip

Theorem 4.6.1 shows that this is a particular case of the Proposition 3.6.1
of [LoMa]. Moreover, the discussion in 3.2 of [LoMa]
establishes a bijection between the even formal solutions to the
equation (5.2) 
and pencils of formal flat connections
$\nabla_0 +\lambda\Cal{A}$ on the trivial vector bundle
with fiber $F$ on the formal completion of 
$T\oplus F$ at zero.

\smallskip

In the remaining part of this section we will establish
certain properties of solutions to (5.2) and 
related structures. As in 3.1--3.2 of [LoMa], we will work
in a more general setting, allowing not necessarily formal manifolds
as our base spaces.

\medskip

{\bf 5.2. Commutativity equations.} Generally, let $M$ be a (super)manifold
in one of the standard categories ($C^{\infty}$, analytic, formal ...).
Denote by $F$ a local system of finite--dimensional
vector (super)spaces on $M$. This is essentially the same as a locally
free sheaf $\Cal{S}$ endowed with a flat connection
$\nabla_0:\,\Cal{S}\to \Omega^1_M\otimes_{\Cal{O}_M}\Cal{S}$:
from $(\Cal{S},\nabla_0)$ one gets $F:=\roman{Ker}\,\nabla_0$,
and from $F$ one gets $\Cal{S}:=\Cal{O}_M\otimes F$ (tensor product
over constants), $\nabla_0 (t\otimes s)=dt\otimes_{\Cal{O}_M}s$
for local sections $t\in \Cal{O}_M,$ $s\in F$.

\smallskip

Here $d$ is the de Rham differential, which extends in the standard
way to $\Omega^*_M$, whereas $\nabla_0$ extends to the whole
tensor algebra of $S$. 

\smallskip

In particular, we have the induced flat connection denoted by the same
letter $\nabla_0$ upon $End\,\Cal{S} =\Cal{O}_M\otimes End\,F.$

\smallskip
 
An even section $B$ of $\Cal{O}_M\otimes\,End\,F$
is called {\it a solution to  the commutativity equations
for $(M,F)$} if 

$$
\nabla_0B\,\wedge\, \nabla_0B=0.
\eqno(5.3)
$$

\smallskip

Our sign conventions
are determined by postulating  that $d$ and $\nabla_0$
are odd. We traditionally denote by $\wedge$ the multiplication
in $\Omega^*_M$ and in the tensor products of $\Omega^*_M$
with other sheaves of algebras, although in the supergeometry
this is slightly misleading: multiplication in $\Omega^*_M$
is supercommutative, not superalternate.

\smallskip

If we choose local coordinates $(t^i)$ in $M$ and a basis of
flat (belonging to $F$) local sections of $\Cal{S}$, $B$
becomes a matrix function $t\mapsto B(t)$,
and $\nabla_0B$ becomes a matrix of 1--forms $dB$ on $M$:
$(\nabla_0B)_i^j=dB_i^j$ so that (5.2) is a particular case of (5.3).

\smallskip

Putting $dB=\sum_i dt^ib_i(t)$ where $t=(t^i)$, so that  
$$
b_i(t):=\frac{\partial B}{\partial t^i}\in \Cal{O}_M \otimes \,\roman{End}\,F
$$
we can easily check that (5.3) can be written in the form of commutativity equations
$$
\forall\,i,j,\ [b_i(t),b_j(t)]=0
\eqno(5.4)
$$
where the brackets denote the supercommutator.
Thus $b_i$ span over $\Cal{O}_M$ a sheaf of (super)abelian
Lie subalgebras of $\Cal{O}_M \otimes \,\roman{End}\,F$. This sheaf
is intrinsically associated with $B$; we will denote it $DB$.

\medskip

\proclaim{\quad 5.2.1. Definition} (i) A solution $B$ to the commutativity
equations as above is called maximal one, if $DB$ is a maximal subsheaf
of (super)abelian
Lie subalgebras of $\Cal{O}_M \otimes \,\roman{End}\,F$ in the following sense:
any local section $c(t)\in \Cal{O}_M\otimes \,\roman{End}\,F$ (super)commuting
with $DB$ belongs to $DB$.

\smallskip

(ii) A maximal solution is called strictly maximal if $DB$ is locally freely generated by 
$b_i(t)$.
\endproclaim

\smallskip

Notice that if $T$ is
the formal completions of a linear space at zero,
maximality of $DB$ is equivalent to the maximality at zero.

\medskip

{\bf 5.2.2. Pullback.} Any morphism $\varphi :\,M^{\prime}\to M$ 
and any solution $B$ to the commutativity equations for $(M,F)$ produces the pullback
solution $\varphi^*(B)$ for $(M^{\prime},\varphi^*(F)).$. Of course, if $M$ and $M^{\prime}$
are formal completions of  linear spaces at zero, $\varphi$
is generally given by formal series with vanishing
constant terms. If $\varphi$ is a closed embedding, we
call $\varphi^*(B)$ the restriction of $B$, and 
conversely, we call $B$ a continuation of $\varphi^*(B)$
(both with respect to $\varphi$). 

\smallskip

In particular, the automorphism group of $(M,F)$ acts upon the space of
solutions. In the formal case,
this is the group of formal invertible coordinate changes.
Hence the linear structure of $M$ plays no role, and we can
simply speak about formal manifolds.

\medskip

\proclaim{\quad 5.2.3. Proposition} Let $B$ be a 
solution to the commutativity equations over a
manifold $M$, and $B^{\prime}$
its continuation to $M^{\prime}$ with respect to
a closed embedding $\iota :\,M\to M^{\prime}$. Denote
by $\widehat{M}^{\prime}$ the formal completion of $M^{\prime}$
along $M$, and by $\widehat{B}^{\prime}$ the solution induced
by $B^{\prime}$ on it. 

\smallskip

(i) If $B$ is maximal, everywhere locally over $M$
there exists a formal projection $\varphi :\, \widehat{M}^{\prime}\to M,\,
\varphi \circ \iota =id_M$, such that $\widehat{B}^{\prime}=\varphi^*(B).$

\smallskip

(ii) If $B$ is strictly maximal, these local formal projections are unique
and hence glue together to a global formal projection
$\hat{M}^{\prime}\to M$ with the same property.
\endproclaim

\smallskip

{\bf Proof.} Working in a local chart, choose local 
coordinates 
$$
(t,\theta )=(t^1,\dots ,t^m; \theta^1,\dots ,\theta^n )
$$ 
on $M^{\prime}$ such that $M$ in $M^{\prime}$
is given by the equations $\theta =0$. Choosing also a basis
of sections of $F$, consider $B$ and $B^{\prime}$ as
matrix functions. 
Denote by $\theta^{(k)}$ (super)symmetric monomials in $\theta$
where $(k)$ are polydegrees. We have
$$
dB\wedge dB=0,\ dB^{\prime}\wedge dB^{\prime}=0,\ B^{\prime}(t,0)=B(t).
\eqno(5.5)
$$
Let $|k|$ be the total degree (sum of coordinates) of the polydegree $(k).$
To prove the statement (i), we have to find a family of
local functions  $(\lambda_{(k)}^i(t))$, $i=1,\dots ,m$, $|k|\ge 1$, on $M$,
such that putting
$$
B^{\prime}_N(t,\theta ):= B(t+\sum_{|k|\le N} \lambda_{(k)}(t) \theta^{(k)}),
\eqno(5.6)
$$
we have
$$
B^{\prime}_N(t,\theta )\, \equiv\, B^{\prime}(t,\theta )\,\roman{mod}\, \theta^{\ge N+1},
\eqno(5.7)
$$
where $\theta^{\ge N+1}$ is the ideal generated by $\theta^{(k)}$ with $|k|\ge N+1.$
In fact, if we find such a family, the local formal projection $\varphi$ we look for
is given by $\varphi^*(t^i)=t^i+\sum_{(k)} \lambda_{(k)}^i(t) \theta^{(k)}$.
Uniqueness of this family means uniqueness of the formal projection.

\smallskip

Obviously, $B^{\prime}_0(t,\theta ) = B(t)$ is the only possible choice.
In order to find $B^{\prime}_1(t,\theta )=B(t+\sum_{j=1}^n \lambda_j(t)\theta^j)$
satisfying (5.7) notice that
$$
B(t+\sum_{j=1}^n \lambda_j(t)\theta^j)\,\equiv
B(t) +\sum_{i=1}^m\sum_{j=1}^n\lambda_j^i(t)\theta^jb_i(t)\ \roman{mod}\, \theta^{\ge 2}
$$
whereas
$$
B^{\prime}(t,\theta ) \equiv B(t) +\sum_{j=1}^n  \frac{\partial B^{\prime}}{\partial \theta^j}
\Bigm|_{\theta =0}
\,\theta^j\
\roman{mod}\, \theta^{\ge 2}.
$$
Hence $\lambda_j^i(t)$ exist iff all
$\dfrac{\partial B^{\prime}}{\partial \theta^j}
\Bigm|_{\theta =0}$
belong to $DB$ that is, (super)commute with $b_i(t).$
But this follows from $dB^{\prime}\wedge dB^{\prime}=0$ evaluated at
$\dfrac{\partial}{\partial t^i}\wedge \dfrac{\partial}{\partial \theta^j}.$ 
Moreover, $\lambda_j^i(t)$ are unique, if $DB$ is freely generated by 
$(b_i(t))$.

\smallskip

This reasoning can be generalized to provide an inductive step
from $N$ to $N+1$.

\smallskip

Namely, assume that we already found $B^{\prime}_N(t,\theta )$
of the form (5.6) satisfying (5.7), $N\ge 1$. We need a vector  
$X_{N+1}=(X_{N+1}^j)$, $j=1,\dots ,n$, whose coordinates are forms
of degree $N+1$ in $\theta$ with coefficients
depending on $t$ such that
$$
B^{\prime}_{N+1}(t,\theta ):=B(t+\sum_{|k|\le N} \lambda_{(k)}(t) \theta^{(k)}+X_{N+1})\equiv
$$
$$
\equiv B^{\prime}_{N}(t,\theta ) + \sum_j X_{N+1}^jb_j(t+\sum_{|k|\le N} \lambda_{(k)}(t) \theta^{(k)})\
\roman{mod}\, \theta^{\ge N+2}
$$
$$
\equiv B^{\prime}_{N}(t,\theta ) + \sum_j X_{N+1}^jb_j(t)\
\roman{mod}\, \theta^{\ge N+2}
$$
satisfies (5.7) with $N$ replaced by $N+1$. Let
$$
B^{\prime}(t,\theta ) \equiv B^{\prime}_N(t,\theta )+ Y_{N+1} \ \roman{mod}\,\theta ^{\ge N+2}
$$ where $Y_{N+1}$ is a form of  degree $N+1$ in $\theta$ with coefficients
depending on $t$. In order to establish the existence of
$X_{N+1}$ it suffices to check as above that the coefficients
of $Y_{N+1}$ commute with $b_j(t)$. In fact, since $dB^{\prime}\wedge dB^{\prime}=0$, we have
$$
d(B^{\prime}_N +Y_{N+1})\wedge d(B^{\prime}_N +Y_{N+1}) \equiv 0\ \roman{mod}\,\theta^{\ge N+1}
$$
and since moreover $dB^{\prime}_N\wedge dB^{\prime}_N=0$ and $N\ge 1$, it follows that
$$
dY_{N+1}\wedge dB^{\prime}_N + dB^{\prime}_N\wedge dY_{N+1} 
\equiv 0\ \roman{mod}\,\theta^{\ge N+1} .
$$
Hence the coefficients of all derivatives $\dfrac{\partial}{\partial \theta^j}\,Y_{N+1}$ belong to $DB$.
From Euler's formula it follows that the same holds for the coefficients
of $Y_{N+1}$. 
Uniqueness in the case of strict maximality follows for this inductive step as well.

\medskip

{\bf 5.2.4. Primitive vectors.} We keep the notation described above. Consider fibers of
$F$ as  (super)manifolds endowed with a linear structure. We can cover $M$ by open
submanifolds such that over any chart $U$ we have a canonical
trivialization $F|_U=F_0\times U$, where $F_0$ is a fiber, or else
the space of sections of $F$ over U.

\smallskip 

An even section $h\in \Gamma (M,F)$ 
is called {\it primitive (for the solution $B$)} if the map of (super)manifolds
$U\to F_0:\, t\mapsto B(t)h$ is a local isomorphism everywhere on $M$. An evident necessary
condition for the existence of a primitive vector is the coincidence of
(super)dimensions of $M$ and fibers of $F$.

\medskip

{\bf 5.3. Oriented associativity equations.} Consider now
a (super)manifold $M$ endowed with an affine flat structure.

\smallskip

By definition, such a structure is given by
the subsheaf $\Cal{T}_M^f\subset \Cal{T}_M$ of flat vector fields
which form a local system of linear spaces {\it and abelian Lie superalgebras}
(with respect to the supercommutator of vector fields) such that
$\Cal{T}_M=\Cal{O}_M\otimes \Cal{T}_M^f.$ Let $\nabla_0:\,\Cal{T}_M\to 
\Omega_M^1\otimes_{\Cal{O}_M}\Cal{T}_M$ be the associated connection.

\smallskip

Consider a vector field $A$ on $M$. Its covariant differential
$\nabla_0A$ belongs to $\Omega^1_M\otimes_{\Cal{O}_M}\Cal{T}_M.$
We have the standard isomorphism
$$
j:\, \Omega^1_M\otimes_{\Cal{O}_M}\Cal{T}_M \to End_{\Cal{O}_M}(\Cal{T}_M).
$$ 
Put $B:=j(\nabla_0A).$ We are now in a position to write
the commutativity equations for $(M,\Cal{T}_M^f)$ and $B$.

\smallskip

A concrete way of fixing a flat structure consists in giving an atlas with affine linear
transition functions between local coordinates of its charts. Such
local coordinates $(t^i)$ are then called flat, and $\Cal{T}_M^f$
is locally generated by the dual vectors $\partial_i=\partial /\partial t^i.$
When we write the field $B=B(t)$ as a matrix, we always use $(\partial_i)$
as a basis of local sections.

\medskip

\proclaim{\quad 5.3.1. Definition} (i) A solution to the oriented associativity equations
is an even vector field $A$ on $M$ with the following property. Write
$A =\sum_c A^c\partial_c$ in a local flat coordinate system and put
$A_{ab}{}^c := \partial_a\partial_b A^c$. Then the composition law
$$
\partial_a\circ \partial_b = \sum_c A_{ab}{}^c\partial_c
\eqno(5.8)
$$
extends to the associative (super)commutative $\Cal{O}_M$--bilinear
multiplication on $\Cal{T}_M$.

\smallskip

(ii) An even flat vector field $\partial_0$ (contained among $\partial_i$) 
is called a flat identity (for $A$),
if $A_{a0}{}^c=\delta_a^c$, or in other words, if $e\circ X=X$
for each vector field $X$. 
\endproclaim

\medskip

\proclaim{\quad 5.3.2. Proposition} Let $A$ be a solution to the 
oriented associativity equations with flat identity $e$. 
Define  $B:=j(\nabla_0A)$ as above and put $h=e$.

\smallskip

Then $B$ is a solution to the commutativity equations with primitive vector $h$.

\smallskip
The map $t\mapsto B(t)e$ establishes
 a local embedding of $M$ into a fiber of $F$ such
that the tautological flat affine structure on $F$ induces the initial affine
structure on $M$. 
\endproclaim
 
\smallskip

{\bf Proof.} In the notations of 5.3.1, we have
$$
\nabla_0A=\sum_{b,c}dt^b\partial_bA^c\otimes \partial_c =
$$
$$
=\sum_{b,c}(-1)^{(|t_b|+1)(|t_b|+|t_c|)} \partial_bA^c dt^b\otimes \partial_c.
$$
where $|t^a|$ denotes the parity of $t^a.$
The last line allows us to compute $j(\nabla_0A)$ as a matrix $B$:
$j(dt^b\otimes\partial_c)$ considered as endomorphism
of $\Cal{T}_M^f$ maps $\partial_a$ to $\delta_a^b\partial_c$ so that
$$
B_b{}^c=(-1)^{(|t_b|+1)(|t_b|+|t_c|)} \partial_bA^c.
\eqno(5.9)
$$
Then a direct calculation shows that
the equation $dB\wedge dB=0$ becomes
$$
\forall\, a,b,c,f,\qquad \sum_e A_{ab}{}^eA_{ec}^f=
(-1)^{|t^a|(t^b|+|t^c|)}\sum_e A_{bc}{}^eA_{ea}{}^f.
$$
These are precisely the associativity equations:
cf. e.~g.  [Ma1], pp. 19--20. Since this formulation gives
a coordinate free description of the operation $\circ$,
it does not depend on the choice of local flat coordinates
and agrees on intersections.
Supercommutativity of $\circ$ follows from the fact
that $\Cal{T}_M^f$ is a superabelian Lie algebra
so that $A_{ab}{}^c$ is symmetric with respect to its subscripts.
\smallskip

Let now $t=(t^0,t^1,\dots ,t^m)$ be a flat local coordinate system
such that $e=\partial_0.$ We have $\partial_a(A_0^c)=\delta_a^c$,
hence $A_0^c=t^c+a^c$ where $a^c$ are constants. Therefore 
$$
B(t)e= \sum_c (-1)^{|t^c|}(t^c+a^c)\partial_c.
\eqno(5.10)
$$ 
This proves the last statement 
of the Proposition.

\medskip

The converse statement also holds, at least locally.
The point is that for a given $B$, if the equation $B=j(\nabla_0A))$
can be solved for $A$ at all, then we can get other
solutions by adding to $A$ any flat vector field,
so there may be an obstruction for finding
a global solution.

\medskip

\proclaim{\quad 5.3.3. Proposition} Let $B$ be a solution to the commutativity
equations for $(M,F)$, and let $h$ be a primitive vector for $B$.
Working locally, induce
a flat affine structure on $M$ from a fibre  $F_0$ with the help of the map
$t\mapsto B(t)h$.

\smallskip

Then in any local flat coordinate system $(t^c)$ we can define
a vector field $\sum_c A^c\partial_c$ such that 
$$
\partial_bA^c = (-1)^{(|t^b|+1)(|t^b|+|t^c|)}B_b{}^c.
$$
(cf. (5.9)).

This field is a solution to the oriented associativity equations
with flat identity $e$ which is the pullback of $h$ (considered
as a tangent vector to a fiber of $F$) with respect to the map $t\mapsto B(t)h$. 
\endproclaim

{\bf Proof.} If we prove the existence of functions $A^c$, the rest will follow
from the proof of the previous Proposition. The equations
for $A^c$ are equivalent to 
$$
dA^c=\sum_b (-1)^{(|t^b|+1)(|t^b|+|t^c|)}dt^bB_b{}^c
$$ 
and their
integrability is equivalent to the closedness of all forms 
$$
\omega^c:=(-1)^{(|t^b|+1)(|t^b|+|t^c|)}\sum_b dt^bB_b{}^c.
$$
We will show that the latter in appropriate coordinates 
is expressed by the equations $dB\wedge d(Bh) =0$
which follow from $dB\wedge dB=0$ because $h$ is flat. 

\smallskip

In fact, choose a basis of flat vector vector fields of $F_0$
containing $h$ and the dual coordinate system $(t^c)$ such that
$h=\partial_0.$ Identifying  $M$ (or its local chart) with a subdomain
of $F$ via $t\mapsto B(t)h$ we can consider $(t^c)$ as flat coordinates on $M$,
and in these coordinates the map $t\mapsto B(t)h$ is given
by the formula (5.10). Hence 
$$
d(B(t)h)=\sum dt^b\otimes \partial_b
$$
and
$$
dB\wedge d(B(h))=\sum_c\left[\sum_b (dB)_b{}^c\wedge dt^b\right]\otimes\partial_c=
\sum_c (-1)^{|t^c|}d\omega^c \otimes\partial_c
$$
This completes the proof.

\medskip

{\bf 5.4. Formal solutions to the commutativity equations corresponding
to $\Cal{L}$--algebras.} Let us now return to the situation described in
the Theorem 5.1.1. The formal series $\Cal{B}$ is a solution to
the commutativity equations for the formal manifold
$M$ which is the completion of $T\oplus F$ at zero, and trivial
local system with the fibre $F$.
 Denote by $\Cal{T}$ (resp. $\Cal{F}$) formal
completions of $T$ (resp. $F$) at zero.
They are embedded as closed formal submanifolds in $M$.

\smallskip

We can restrict $\Cal{B}$ to $\Cal{T}$. Assume that this restriction $\Cal{B}_T$ is maximal.
Then in view of the Proposition 5.2.3, $\Cal{B}$ is a pullback of  $\Cal{B}_T$ with respect to
a formal projection $\varphi :\,M\to \Cal{T}$. 

\smallskip

We can also restrict $\Cal{B}$ to $\Cal{F}$.
Clearly, $\Cal{F}$ is endowed with a formal flat structure coming from  $F$.
The restriction $\Cal{B}_F$ of $\Cal{B}$ to $\Cal{F}$  produces a solution to
the oriented associativity equations with base $\Cal{F}$. If the latter admits a primitive
vector $h$, we will say that $h$ is weakly primitive for $\Cal{B}$.
\medskip

\proclaim{\quad 5.4.1. Proposition} Fix two linear superspaces $T$ and $F$
and the following additional data:

\smallskip

(i) A maximal solution $\Cal{B}_1$ to the commutativity equations with the base
$\Cal{T}$ and fiber $F$.

\smallskip

(ii) A solution $\Cal{B}_2$ to the commutativity equations with the base
$\Cal{F}$ and fiber $F$ which comes from a solution to the oriented
associativity equations on $\Cal{F}$ with flat structure induced by $F$,
admitting a flat identity $e$.

\smallskip

In this case there exists a pair $(\Cal{B},h)$, where 
$\Cal{B}$ is a solution to
the commutativity equations on the formal completion of $T\oplus F$ at zero,
with the fiber $F$, and $h$ is a weakly primitive vector corresponding to $e$,
such that $\Cal{B}_1=\Cal{B}_T$, $\Cal{B}_2=\Cal{B}_F.$ 
This pair is unique.
\endproclaim

\medskip

{\bf Proof.} Let $(t)$ denote some coordinates on $\Cal{T}$ and $(\theta)$
flat coordinates on $\Cal{F}$. The explicit pullback formula
$$
\Cal{B}(t,\theta )=\Cal{B}_2(\theta +\Cal{B}_1(t)h)
$$
produces a solution with necessary properties. It is unique because
of maximality of $\Cal{B}_1$.

\medskip

{\bf 5.5. Compatibility of two tensor products.} Given two structures $I_i$
of (oriented) $\Cal{L}$--algebras upon  $(T_i,F_i)$, $i=1,2$, we can
form the tensor product structure $I=I_1\otimes I_2$ upon
$(T,F)=(T_1\otimes T_2,F_1\otimes F_2)$ imitating the definition 4.2.1
in the nonoriented (cyclic) case. This operation induces
the tensor product $*_{COMM}$ on the formal solutions to the
commutativity equations (COMM). 

\smallskip

There is also an oriented version $*_{ASS}$ of the tensor product
of the formal solutions to the associativity equations ASS (for the cyclic case,
see [Ma1], p. 100).

\smallskip

Considering only solutions to ASS with flat identities, and the associated
solutions to COMM, we may conjecture that the two operations are
compatible. This question was raised in [LoPo] where it was
checked that the answer is positive in several first orders.
\smallskip

Proposition 5.4.1 suggests a strategy for proving
this conjecture at least for the case when the solutions to COMM
are maximal. However, the most natural and general approach is furnished by
the equivalence theorem 5.1.1 and the language of $\Cal{L}$--algebras where
the tensor product is simply induced
by the product in the cohomology of moduli spaces,
resp. coproduct in their homology. We hope to spell out
the details elsewhere.

\bigskip

\centerline{\bf References}

\medskip

[BeMa] K.~Behrend, Yu.~Manin. {\it Stacks of stable maps and Gromov--Witten invariants.} 
Duke Math. Journ., 85:1 (1996), 1--60.

\smallskip

[GeK1] E.~Getzler, M.~Kapranov. {\it Cyclic operads and cyclic homology.}
In: Geometry, Topology, and Physics for Raoul, ed.~by B.~Mazur,
Internat.~Press, Cambridge, MA, 1995, 167--201.

\smallskip

[GeK2] E.~Getzler, M.~Kapranov. {\it Modular operads.}
Compositio Math., 110 (1998), 65--126.

\smallskip

[H] B.~Hassett. {\it Moduli spaces of weighted pointed stable curves.}
Preprint \newline math.AG/0205009

\smallskip

[Ke] S.~Keel. {\it Intersection theory of moduli space
of stable $N$--pointed curves of genus zero}.
Trans.~AMS, 330:2 (1992), 545--574.

\smallskip

[Kn] Knudsen F.F. {\it The projectivity of the moduli space
    of stable curves II.
The stacks $\overline M_{0,n}$}. Math. Scand. 52 (1983), 163--199.

\smallskip

[KoMa] M.~Kontsevich, Yu.~Manin. {\it Gromov--Witten classes, quantum cohomology, and enumerative 
geometry.} Comm. Math. Phys.,
164:3 (1994), 525--562.

\smallskip 

[KoMaK] M.~Kontsevich, Yu.~Manin (with Appendix by
R.~Kaufmann). {\it Quantum cohomology of a product.} Inv. Math., 124, f. 1--3 (Remmert's Festschrift) (1996), 313--339.

\smallskip

[Lo1]  A.Losev,{\it   On ``Hodge'' Topological Strings at genus zero.}
  Pis'ma v ZhETF  65, 374-379 (1997).

\smallskip

[Lo2] A.Losev,  {\it ``Hodge strings'' and elements of  K.~Saito's theory of
the
primitive form.}
In: Proceedings of
Taniguchi Symposium on Topological Field Theory, Primitive Forms
and Related Topics, Kyoto,
Japan, 9-13 Dec 1996, Springer,1998.
Preprint hep-th/9801179

\smallskip

[LoMa] A.~Losev, Yu.~Manin. {\it New moduli spaces of pointed 
curves and pencils of flat connections.}
Michigan Journ. of Math., vol. 48
(Fulton's Festschrift), 2000, 443--472. Preprint math.AG/0001003

\smallskip

[LoPo] A.~Losev, I.~Polyubin. {\it On compatibility of tensor products
on solutions to commutativity and WDVV equations.} JETP Letters,
73:2 (2001), 53--58.

\smallskip

[Ma1] Yu.~Manin. {\it Frobenius manifolds, quantum cohomology,
and moduli spaces.}  AMS Colloquium Publications, vol. 47, Providence, RI, 1999, xiii+303 pp.

\smallskip

[Ma2] Yu.~Manin. {\it Moduli stacks $\overline{L}_{g,S}$.} Preprint
math.AG/0206123

\smallskip

[R] M.~A.~Readdy. {\it The Yuri Manin ring and its
$\Cal{B}_n$--analogue.} Adv. in Appl. Math., 2001.

\enddocument